\documentclass[AMA,Times1COL,fleqn]{WileyNJDv5}

\usepackage{longtable}
\usepackage{amsmath}
\usepackage{amsfonts}
\usepackage{amsthm}
\usepackage{mathtools}

\usepackage{subfig}
\usepackage{lineno}
\usepackage{textcomp}
\usepackage{mathrsfs}
\usepackage{booktabs}

\usepackage{multicol}
\usepackage{xcolor}

\usepackage{tabularx}

\usepackage{amssymb}
\usepackage{physics}
\usepackage{MnSymbol}
\usepackage{empheq}

\usepackage{float}
\usepackage{array,multirow}
\usepackage{color}
\usepackage{tikz}

\usepackage{pgfplots}
\pgfplotsset{compat=1.16}

\definecolor{lightgray}{gray}{0.80}

\definecolor{grey1}{rgb}{0.5, 0.5, 0.5}
\definecolor{green1}{rgb}{0.4660, 0.6740, 0.1880} 
\definecolor{blue1}{rgb}{0, 0.4470, 0.7410} 
\definecolor{red1}{rgb}{0.8500, 0.3250, 0.0980}
\definecolor{yellow1}{rgb}{0.9290, 0.6940, 0.1250}
\definecolor{purple1}{rgb}{0.4940, 0.1840, 0.5560}
\definecolor{lightblue1}{rgb}{0.3010, 0.7450, 0.9330}
\definecolor{bordeaux1}{rgb}{0.6350, 0.0780, 0.1840}
\definecolor{brown1}{rgb}{0.65, 0.16, 0.16}
\definecolor{pink1}{rgb}{1.0, 0.08, 0.58}
\definecolor{orange1}{rgb}{1, 0.65, 0.0}

\definecolor{green2}{rgb}{0.0, 0.4980, 0.0}
\definecolor{burntorange}{rgb}{0.74902,0.341176,0}

\newcommand{\PreserveBackslash}[1]{\let\temp=\\#1\let\\=\temp}
\newcolumntype{C}[1]{>{\PreserveBackslash\centering}p{#1}}
\newcolumntype{R}[1]{>{\PreserveBackslash\raggedleft}p{#1}}
\newcolumntype{L}[1]{>{\PreserveBackslash\raggedright}p{#1}}

\usepackage{color, soul}
\definecolor{blue2}{RGB}{125, 249, 255}

\DeclareRobustCommand{\reviewerI}[1]{{\sethlcolor{pink}\hl{#1}}}
\DeclareRobustCommand{\reviewerII}[1]{{\sethlcolor{yellow}\hl{#1}}}

\DeclareRobustCommand{\changed}[1]{{\sethlcolor{blue2}\hl{#1}}}

\soulregister\reviewerI{1}
\soulregister\reviewerII{1}
\soulregister\changed{1}

\soulregister\ref7
\soulregister\pageref7
\soulregister\eqref7
\soulregister\cite7

\makeatletter
\renewcommand*\env@matrix[1][*\c@MaxMatrixCols c]{%
  \hskip -\arraycolsep
  \let\@ifnextchar\new@ifnextchar
  \array{#1}}
\makeatother

\theoremstyle{plain}

\theoremstyle{definition}

\newcommand{\vect}[1]{\boldsymbol{#1}} 									
\newcommand{\mat}[1]{\mathbf{#1}} 											

\newcommand\abss[1]{\left\vert#1\right\vert}       

\newcommand{\m}{m}																	

\newcommand{\bspline}{N}															

\newcommand{\phic}{\boldsymbol{\varphi}}

\newcommand{\qhat}{q}

\newcommand{\Pd}{\mathcal{P}_{\vect{d}}}
\newcommand{\Hd}{\mathcal{H}_{\vect{d}}}

\newcommand{\iga}{IGA}
\newcommand{\nodal}{Nodal $\mathbb{R}^3$}
\newcommand{\nodalSaddle}{Nodal SPP}
\newcommand{\nodalSaddleRed}{Nodal SPP-reduced}
\newcommand{\nodalPenalty}{Nodal-penalty}

\def\genbox#1#2#3#4#5#6{
    \leavevmode\raise#4bp\hbox to#5bp{\vrule height#5bp depth0bp width0bp
    \pdfliteral{q .5 w \csname #2COLOR\endcsname\space RG
                       \csname #3PDF\endcsname{#5}{#6} S Q
             \ifx1#1 q \csname #2COLOR\endcsname\space rg 
                       \csname #3PDF\endcsname{#5}{#6} f Q\fi}\hss}}

\articletype{Research Article (accepted version)}

\received{Date Month Year}
\revised{Date Month Year}
\accepted{Date Month Year}

\journal{International Journal for Numerical Methods in Engineering}
\copyyear{2025}
\startpage{1}

\graphicspath{{./fig/}}

\begin{document}

\title{A study on nodal and isogeometric formulations for nonlinear dynamics of shear- and torsion-free rods}

\author[1]{Thi-Hoa Nguyen}

\author[1]{Bruno A. Roccia}

\author[2]{Dominik Schillinger}

\author[1]{Cristian G. Gebhardt}

\authormark{Nguyen \textsc{et al.}}
\titlemark{A study on nodal and isogeometric formulations for nonlinear dynamics of shear- and torsion-free rods}

\address[1]{\orgdiv{Geophysical Institute and Bergen Offshore Wind Centre}, \orgname{University of Bergen}, \orgaddress{ \country{Norway}}}

\address[2]{\orgdiv{Institute for Mechanics, Computational Mechanics Group}, \orgname{Technical University of Darmstadt}, \orgaddress{ \country{Germany}}}

\corres{\email{hoa.nguyen@uib.no}}

\abstract[Abstract]{In this work, we compare the nodal and isogeometric spatial discretization schemes for the nonlinear formulation of shear- and torsion-free rods introduced in \cite{gebhardt_2021_beam}. 
We investigate the resulting discrete solution space, the accuracy, and the computational cost of these spatial discretization schemes. 
To fulfill the required $C^1$ continuity of the rod formulation, the nodal scheme discretizes the rod in terms of its nodal positions and directors using cubic Hermite splines. 
Isogeometric discretizations naturally fulfill this with smooth spline basis functions and discretize the rod only in terms of the positions of the control points \cite{nguyen_rod_2024}, which leads to a discrete solution in multiple copies of the Euclidean space $\mathbb{R}^3$. They enable the employment of basis functions of one degree lower, i.e. quadratic $C^1$ splines, and possibly reduce the number of degrees of freedom. 
When using the nodal scheme, since the defined director field is in the unit sphere $S^2$, preserving this for the nodal director variable field requires an additional constraint of unit nodal directors. This leads to a discrete solution in multiple copies of the manifold $\mathbb{R}^3 \times S^2$, however, results in zero nodal axial stress values. 
Allowing arbitrary length for the nodal directors, i.e. a nodal director field in $\mathbb{R}^3$ instead of $S^2$ as within discrete rod elements, 
eliminates the constrained nodal axial stresses and leads to a discrete solution in multiple copies of $\mathbb{R}^3$. 
To enforce the unit nodal director constraint, we discuss two approaches using the Lagrange multiplier and penalty methods. 
We compare the resulting semi-discrete formulations and the computational cost of these discretization variants. 
We numerically demonstrate our findings via examples of a planar roll-up, a catenary, and a mooring line. 
}


\maketitle

\section{Introduction}
  
Nonlinear rods find their broad applications in various areas of science and engineering, 
such as in the analysis of DNA molecules \cite{Schlick1995,Shi1994}, 
microstructures \cite{Cyron_beam2012,Muller2015}, 
the dynamics of cables \cite{Boyer2011,Coyne1990}, 
the mechanical analysis of M\"obius bands \cite{Moore2019}, or the 
stability of elastic knots \cite{Audoly2007,Ivey1999}, among others. 
Particularly in offshore and coastal engineering, cables or rods that are employed for towing or mooring of floating structures, belong to the most important research topics. 
For such applications of long and slender cables, where shear deformations can be neglected, 
shear-free rod models have been shown to 
accurately capture their deformations and behavior, 
and 
hence 
belong to the most relevant and interesting models. 
They are based on the assumption of cross-sections that remain flat and perpendicular to the tangent vector of the rod axis (see, e.g., assumptions of rod formulations studied in\cite{Giusteri2018,OReilly2017}), which 
simplifies and reduces the number of kinematic degrees of freedom, compared e.g. to models capturing shear and/or torsion deformations. 
In general, shear-free models are developed on the basis of the Kirchhoff-Love theory. 
We note that another category of rod models that capture the shear deformations is generally based on the Simo-Reissner theory. 
For a comprehensive overview and review of different rod formulations\footnote{Note that in literature, ``rod'' models/structures are also refered to as ``beams''.} developed based on these two theories, we refer to \cite{Meier2019} and reference therein.

For linear shear-free rods, well-established models are the Euler-Bernoulli and Rayleigh models \cite{Antman1972,HAN1999}. 
For nonlinear shear-free rods, one of the most widely employed models is the Kirchhoff rod model, which can be considered a generalization of the Rayleigh rod model \cite{Antman1974,Antman2005,Langer1996}. 
Various attempts have been made to develop shear-free nonlinear rod formulations and solve the resulting governing equations using finite element methods. 
There exists a variety of rod (beam) elements that are geometrically exact and able to capture large deformations \cite{Boyer2011,Raknes2013,Greco2014,Greco2016,Boyer2004,Bauer2016,weiss2002,Zhang2015}. 
The developed rod elements of these formulations are Kirchhoff rod elements, for which Armero and Valverde give a general historic overview in \cite{Armero2012,Armero2012a}. 
An alternative to the Kirchhoff rod elements is the so-called corotational shear-free beams \cite{Hsiao1994,Hsiao1999,Le2014}. 
We note that in general, the above mentioned rod formulations and elements are based on constrained variational statements, particularly in dynamics problems. 
Due to the non-integrable nature of vanishing shear deformations, 
it is generally not possible to formulate the governing equation of shear-free rods through a truly unconstrained variational statement \cite{Giusteri2018,OReilly2017}. 
In the last few years, there have been attempts to tackle this. 
In \cite{Romero2020}, Romero and Gebhardt developed an unconstrained variational formulation for shear-free Kirchhoff rods but relied on certain simplification hypotheses. 
In their following work \cite{gebhardt_2021_beam}, the authors introduced 
an unconstrained structural model for nonlinear initially straight rods that do not exhibit shear and torsion, where the isotropy of the cross-section no longer plays an essential role. 
This can be considered as a special case of the static and dynamic variational principles for the Kirchhoff rod model developed in \cite{Romero2020}, 
as well as a non-shearable counterpart of the torsion-free beam model introduced in \cite{Romero2014}. 
In this work, we consider the shear- and torsion-free nonlinear rod formulation developed in \cite{gebhardt_2021_beam} that is a Kirchhoff rod model and 
in its simplest representation and does not include any non-integrable constraints such as the one enforcing non-twisting conditions. 
We note that in the context of Kirchhoff rods, the existence of membrane locking \cite{Stolarski1982} is reported, see e.g. \cite{Meier2019,Armero2012,Armero2012a}, 
which can be eliminated using, for instance, the approach of reduced/selective integration (see e.g. \cite{Noor1981,Adam2015,Zou_quadrature2021}) or approaches based on Hu-Washizu or Hellinger–Reissner variational principles (see e.g. \cite{Cannarozzi2008,Choit1995,Lee1993}). 
For an overview of different locking-preventing techniques, we refer to e.g. \cite{Meier2019,Nguyen2022} and references therein.

To solve the governing equations of the rod formulation \cite{gebhardt_2021_beam}, nodal finite elements and isogeometric discretizations have been employed in \cite{gebhardt_2021_beam} and \cite{nguyen_rod_2024}, respectively. 
The rod formulation \cite{gebhardt_2021_beam} requires at least $C^1$ continuity, which is naturally fulfilled using smooth spline basis functions when using isogeometric discretizations \cite{hughes_isogeometric_2005,Cottrell:09.1}. 
This type of discretization has been widely employed for different structures such as rods and beams, e.g. in \cite{Greco2016,Cazzani2016,Bouclier2012,weeger2013isogeometric}, and 
shells e.g. in \cite{Alaydin2021,Benson_shell_2013,Borkovic2022,Echter_shell_2013,Kiendl_shell_2009,Oesterle2022}. 
In \cite{gebhardt_2021_beam}, to achieve $C^1$ continuity, the authors employed nodal finite elements and discretized the rod configuration in terms of its nodal positions and nodal directors using standard cubic Hermite splines. 
The obtained discrete solution belongs to multiple copies of the manifold $\mathbb{R}^3 \times S^2$ since the nodal directors are constrained to the unit sphere $S^2$, i.e. unit nodal directors. 
We note that nodal finite elements are also employed for various structures including 
rods and beams \cite{Boyer2004,Meier2014,Zhao2012,Kim1998}, plates 
and shells \cite{Romero2002,Kloppel_scaled_director2011,Kim1996,Bischoff2001,Bathe1985}.

In this work, we give an overview and attempt to gain a deeper understanding of the discretization schemes based on nodal and isogeometric finite elements for the rod formulation \cite{gebhardt_2021_beam}. 
We discuss and show that when using the nodal discretization scheme, 
since the defined director field of the continuous rod configuration lives in the unit sphere $S^2$, 
preserving this for the nodal director variable field requires an additional constraint of unit nodal directors and leads to zero nodal axial stress values. 
We discuss two approaches to enforce the unit nodal director constraint:  
the Lagrange multiplier method that is also employed in \cite{gebhardt_2021_beam} ((see, e.g., also \cite{Fernandez-Mendez2004,Codina2015})), and  
the penalty method (see e.g. \cite{Pasch2021,Lu2019,Schillinger2016}). 
While the former leads to a discrete solution in multiple copies of the manifold $\mathbb{R}^3 \times S^2$, the latter strictly leads to a solution in multiple copies of the Euclidean space $\mathbb{R}^3$ that is the same when using isogeometric discretizations \cite{nguyen_rod_2024}. 
Allowing arbitrary length for the nodal directors, i.e. nodal directors in $\mathbb{R}^3$ instead of $S^2$ as the director field within rod elements, eliminates constrained nodal axial stress values and leads to a discrete solution in multiple copies of $\mathbb{R}^3$. 
Moreover, we highlight that using isogeometric discretizations enables the employment of basis functions of one degree lower, i.e. quadratic $C^1$ splines, and possibly reduces the number of degrees of freedom. 
We discuss the resulting semi-discrete formulation and matrix equations of each discretization variant. 
Enforcing the unit nodal director constraint using the Lagrange multiplier method leads to matrix equations in the form of a saddle-point problem, where the unknown Lagrange multipliers can be eliminated using nullspace method. 
We discuss and compare the computational cost required for each variant. 
We numerically illustrate via the convergence study of a planar roll-up that preserving the nodal directors in the unit sphere leads to better accuracy in the deformations in different error norms. 
Our results of this pure bending example imply the effect of membrane locking on the stress resultants obtained with any of the studied discretizations. 
Via a static and dynamic example of cables commonly employed as mooring lines, 
we illustrate that all formulations approximately lead to the same final deformed configuration. 
For the static example, they lead to the same stress resultants, except axial stress resultants with zero nodal values when enforcing the unit nodal director constraint using the Lagrange multiplier method. 
For the dynamic example, cubic $C^1$ isogeometric discretization leads to bending moments with larger oscillations and slightly larger responses, which can be due to remaining outliers and/or high-frequency modes (see also discussions in \cite{nguyen_rod_2024}). 
Via these examples, we also numerically illustrate the computational cost required for each formulation in terms of the maximum number of iterations and averaged computing time per iteration. 
We show that  
on the one hand, all formulations require generally the same number of iterations, on the other hand, 
cubic $C^1$ isogeometric discretization requires the least time per iteration, with or without outlier removal. 
Using any of the formulations based on the nodal scheme requires approximately the same computing time, except the reduced saddle-point problem 
which requires significantly more time on fine meshes due to matrix reassembly in each iteration. 
We note that for our computations, we employ the same implicit time integration scheme as in \cite{gebhardt_2021_beam,nguyen_rod_2024}, which is a hybrid combination of the midpoint and trapezoidal rules.

The structure of the paper is as follows: 
In Section \ref{sec:preliminaries}, we briefly review the considered rod formulation introduced in \cite{gebhardt_2021_beam} and the employed implicit time integration scheme for our computations in Section \ref{sec:results}. 
In Section \ref{sec:discretizations}, we discuss possible discrete solution spaces when using the isogeometric and nodal discretization schemes. 
We also give an overview of the corresponding semi-discrete formulations for each discretization variant, including those with two different approaches enforcing the unit nodal director constraint. 
In Section \ref{sec:matrix_form}, we discuss and compare the resulting matrix equations of these variants, for which we discuss the computational cost in Section \ref{sec:computational_cost}.
In Section \ref{sec:results}, we numerically illustrate our findings via examples of a planar roll-up, a catenary, and a mooring line. 
In Section \ref{sec:summary}, we summarize our results and draw conclusions.

\section{Preliminaries}\label{sec:preliminaries}

In this section, we briefly review the variational formulation of nonlinear shear- and torsion-free rods in a continuous setting introduced in \cite{gebhardt_2021_beam}. 
We then briefly recall and discuss the most important properties of the employed implicit time integration scheme for our computations in this work.

\subsection{Variational formulation}\label{sec:pre_var_form}

Let the curve $\phic$ now be the configuration of Kirchhoff rods, dependent on the arc-length $s$ and time $t$, $\phic\,=\, \phic(s,\,t)$, $(s,\,t) \,\in\, [0,\,L] \, \times \, [0,\,T]$, that are initially straight, shear-, torsion-free, and transversely isotropic \cite{gebhardt_2021_beam}. 
Next, let us consider the following set for the rod configurations: 
\begin{equation}\label{eq-manifold}
    \mathcal{D} := \left\{ \phic \, \in \, C^2 \left( [0, \,L], \mathbb{R}^3 \right), \;
    \abss{\phic^\prime} \, > \, 0, 
    \phic(0,\,t) = \mathbf{0}, \;
    \phic^\prime \, (0,\,t) \, = \, \vect{E}_3 \right\} \, ,
\end{equation}
where $C^2[0, \,L]$ is the space of $C^2$ continuous functions on $[0, \,L]$, 
$\vect{E}_i$, $i=1,2,3$, are the canonical Cartesian basis of $\mathbb{R}^3$. 
For simplicity and concreteness, 
we adopt for the notation here the clamped boundary condition at $s=0$.

We recall, from \cite{gebhardt_2021_beam}, 
the strong form of the equations of motion governing the space-time evolution for the Kirchhoff rod:
\begin{equation}\label{s-eom}
    \vect{n}^\prime + \left(\, \frac{1}{\abss{\phic^\prime}} \, \vect{d} \,\times\, \nabla_{\vect{d}^\prime} \, \vect{m} \, \right)^\prime \, = \, 
    A_\rho\,\ddot{\phic} \,+\, \left(\, \frac{1}{\abss{\phic^\prime}} \, \vect{d} \,\times\, I_\rho \, \nabla_{\dot{\vect{d}}} \, \dot{\vect{d}} \right)^\prime \,-\, \vect{f}^{\text{ext}} \,,
\end{equation}

\noindent
where $\vect{n}$ and $\vect{m}$ are the stress measures, defined as: 
\begin{align}\label{eq-stress}
    \vect{n} \,=\, EA \, \vect{\epsilon} \,, \qquad \vect{m} \,=\, EI\, \vect{\kappa} \,,
\end{align}

\noindent
respectively, which are conjugated with the following strain measures:
\begin{align}\label{eq-strain}
    \vect{\epsilon} \,:=\, \phic^\prime \,-\, \vect{d} \,, \qquad \vect{\kappa} \,:=\, \vect{d} \,\times\,\vect{d}^\prime \,.
\end{align}

\noindent
Here, 
$E$, $A$, $I$ are the Young's modulus, cross-sectional area, and the moment of inertia of the rod, respectively, and $\vect{d}$ is the director of the curve $\phic$: 
\begin{align}\label{eq-triad}
\hspace{-0.2cm}
    \vect{d} \, := \, \frac{\phic^\prime}{\abss{\phic^\prime}} \, , \quad 
\end{align} 
which is well-defined everywhere along $\phic$. 
The director $\vect{d}$ lives in the unit sphere $S^2 \,:=\left\{ \, \vect{d} \,\in\, \mathbb{R}^3 \, \vert \right.$ $\left. \vect{d} \,\cdot\, \vect{d} \,=\, 1 \, \right\}$ that is a nonlinear, smooth, compact, two-dimensional manifold with no group structure \cite{Eisenberg1979APO,Romero2017}. 
The tangent bundle associated with $S^2$ is also a manifold, which is given by
$T S^2 \,:=$ 
$\left\{ \, (\vect{d},\vect{c}) \,\in\, S^2 \,\times\,\mathbb{R}^3 \,, \vect{d} \,\cdot\, \vect{c} \,=\, 0 \, \right\}$.

Here, $A_\rho$ and $I_\rho$ are the mass per unit length and the inertia density, respectively, i.e. $A_\rho \,=\, \rho\,A$ and $I_\rho \,=\, \rho\,I$, where $\rho$ is the mass density, $A$ the cross-section area and $I$ the moment of inertia of the rod. 
$\vect{f}^{\text{ext}}$ is the external generalized forces, and  
the dot notation in the superscript 
denotes the derivative with respect to time $t$, i.e. $\dot{(\cdot)} \,=\, \partial (\cdot) / \partial \, t$. 
We note that since the director $\vect{d}$ is well-defined along the rod $\phic \,\in\, \mathcal{D}$ (see also \eqref{eq-manifold}), the strain measures \eqref{eq-strain} are also well-defined at every point of the rod.

At time $t=0$, we require the following initial conditions: \\
\begin{displaymath}
    \begin{aligned}
        & \phic \,=\, \phic_0 \quad & \text{on } (s,\,t) \,\in\, [0,\,L] \, \times \, {0} \,, \\ 
        & \dot{\phic} \,=\, \mat{v}_0 \quad & \text{on } (s,\,t) \,\in\, [0,\,L] \, \times \, {0} \,.
    \end{aligned}
\end{displaymath}

Additionally, we require at all times the following boundary conditions; for instance, clamped-free ends:
\begin{subequations}
    \begin{align}
        \text{on } & (s,\,t) \,\in\, {0} \, \times \, [0,\,T]: \,  \phic \,=\, \mat{0}\,, \qquad 
        \phic^\prime \,=\, \vect{E}_3 \,, \\
        \text{on } & (s,\,t) \,\in\, {L} \, \times \, [0,\,T]: \, \vect{n} \,+\, \frac{1}{\abss{\phic^\prime}} \, \vect{d} \,\times\, \left(\, \nabla_{\vect{d}^\prime} \, \vect{m} \,-\, I_\rho \, \nabla_{\dot{\vect{d}}} \, \dot{\vect{d}} \,\right) \,=\, \mat{0} \,, \\
        & \hspace{3.8cm} \frac{1}{\abss{\phic^\prime}} \, \vect{d} \,\times\, \vect{m} \,=\, \mat{0} \,.
    \end{align}
\end{subequations}

According to \cite{gebhardt_2021_beam}, the weak form corresponding to \eqref{s-eom} is then: 
\begin{equation}\label{w-eom}
    \begin{split}
        \int_0^S \, & \delta\phic \,\cdot\, \left( 
            \mathcal{M}\left(\phic^\prime\right)\,\hat{\nabla}_{\dot{\phic}} \, \dot{\phic} \,+\,  \mathcal{B}\left(\phic^\prime,\,\phic^{\prime\prime}\right)^T \, \vect{\sigma} \,-\, \vect{f}^{\text{ext}} \,\right) \, \mathrm{d} \, s \,=\, 0\,,
    \end{split}    
\end{equation}

\noindent
where the mass operator, $\mathcal{M}$, and the linearized strain operator, $\mathcal{B}$, are given by:
\begin{align}
    & \mathcal{M} \,=\, \mathcal{M}\left(\phic^\prime\right) \, :=\, A_\rho \, \mat{I} \,+\,  (\cdot)^{\prime\,T} \, I_\rho \, \frac{1}{\abss{\phic^\prime}^2} \, \Pd \, (\cdot)^\prime \, \label{mass-op}\\
    & \mathcal{B} \,=\, \mathcal{B}\left(\phic^\prime,\,\phic^{\prime\prime}\right) \,:=\, 
    \begin{bmatrix}
        \mat{I} \,-\, \frac{1}{\abss{\phic^\prime}} \, \Pd  & \mat{0} \\
        -\frac{1}{\abss{\phic^\prime}^2} \, \left[\phic^{\prime\prime}\right]_\times \, \Hd     & \frac{1}{\abss{\phic^\prime}} \, \left[\vect{d}\right]_\times 
    \end{bmatrix} \; \begin{bmatrix}
        (\cdot)^\prime \\ (\cdot)^{\prime\prime}
    \end{bmatrix} \,. \label{B-op}
\end{align}
Here, 
$\boldsymbol{\sigma} \,:=\, [\vect{n} \quad \vect{m}]^T$, 
$\Pd$ is the orthogonal projection operator, $\Pd\,:=\, \mat{I} \,-\, \vect{d} \,\otimes\, \vect{d}$, 
$\Hd$ is the Householder operator, 
$\Hd \,:=\, \mat{I} \,-\, 2 \, \vect{d} \,\otimes\, \vect{d}$, and $[\vect{a}]_\times$ denotes the skew-symmetric matrix of a vector $\vect{a} \,=\, \left[a_1 \quad a_2 \quad a_3\right]^T$, i.e.:
\begin{align}
    [\vect{a}]_\times \,=\, \begin{bmatrix}
        0 & -a_3 & a_2 \\ a_3 & 0 & -a_1 \\ -a_2 & a_1 & 0
    \end{bmatrix} \,. \nonumber
\end{align}
The field covariant derivative $\hat{\nabla}_{(\cdot)} \, ( \cdot)$ is the extension of the covariant derivative, acting on smooth fields.

\subsection{Time integration scheme}\label{sec:pre_time_integration}

For our computations in Section \ref{sec:results}, 
we apply the same implicit scheme that is employed in \cite{gebhardt_2021_beam,nguyen_rod_2024}, which is  
a hybrid combination of the midpoint and trapezoidal rules. 
Such an implicit scheme 
achieves second-order accuracy, approximately preserves energy, and exactly preserves the linear and angular momentum \cite{gebhardt_implicit_2020,guo_time_int_2022,Wen_time_int_2022}. 
We note that to exactly preserve the energy the strain measures must be quadratic \cite{Gonzalez_timeInt1996}. 
The employed scheme is based on the one introduced in \cite{gebhardt_implicit_2020} which is closely related to 
the Energy-Dissipative-Momentum-Conserving method \cite{Romero2002,Armero2001,Armero2001a,Armero2003}. 
In this work, to serve our objective of comparing the nodal and isogeometric discretization schemes, we eliminate the dissipation terms of the original scheme \cite{gebhardt_implicit_2020}. This enables us to observe the occurrence and investigate the effects of all contents of the response, including the spurious high-frequency contents (see also \cite{nguyen_rod_2024}). 
For more details on how the chosen implicit time integration scheme is applied to the semi-discrete rod formulation, we refer the readers to \cite{gebhardt_2021_beam,nguyen_rod_2024}.

\section{Nodal and isogeometric spatial discretizations}\label{sec:discretizations}
  
In this section, we discuss and compare the spatial discretization of the nonlinear rod formulation reviewed in the previous section, using either the isogeometric or nodal finite elements. 
We highlight the different resulting discrete solution spaces: using isogeometric discretizations leads to a discrete solution in multiple copies of the Euclidean space $\mathbb{R}^3$. 
Using nodal discretizations leads to a discrete solution in multiple copies either of $\mathbb{R}^3$ or the manifold $\mathbb{R}^3 \times S^2$, where $S^2$ denotes the unit sphere. 
While the latter preserves the director field in $S^2$  
along the complete rod, including the nodes, however, leads to zero nodal axial stress, the former leads to a discrete solution in multiple copies of $\mathbb{R}^3$ that is the same as using isogeometric discretizations and allows non-zero nodal axial stress.
We discuss the unit nodal director constraint and its enforcement to preserve  
the nodal director field in $S^2$, 
as well as the resulting semi-discrete formulation.  
We start with the spatial discretization of the continuous rod configuration.

\subsection{Spatially discrete rod configuration}\label{sec:discrete_rod_space}

The rod formulation \eqref{w-eom} requires discretizations of at least $C^1$-continuity. 
As discussed in \cite{nguyen_rod_2024}, this can be naturally fulfilled using isogeometric discretizations based on smooth spline basis functions, blue{$\bspline_i$}, $1 \leq i \leq \m$, of degree $p$ and conitnuity $C^r$, $1 \leq r \leq p-1$. Here, $\m$ denotes the number of the spline basis functions spanning the corresponding basis. 
The rod configuration $\phic(s,t) \in \mathcal{D}$ \eqref{eq-manifold} can be spatially discretized as follows:
\begin{align}\label{eq-discretize_iga}
    \phic(s,\,t) \,\approx\, \phic_h \,(s,\,t) \,=\, 
    \sum_i^{\m} \, \bspline_i \, (s) \, \vect{x}_i\,(t) \,=\, \mat{\bspline}(s) \, \vect{\qhat} \,, 
\end{align}
where 
$\phic_h\,=\,\phic_h\,(s,\,t) \in \mathbb{R}^3$ denotes the discrete rod configuration in space, which is expressed in terms of the time-dependent position of the $i^{\text{th}}$ control point, 
$\vect{x}_i \, \in \mathbb{R}^3$, and 
$\vect{\qhat}\,=\,\vect{\qhat}(t) \in (\mathbb{R}^3)^{\m}$ is the vector of unknown time-dependent coefficients. 
In this work, we refer to $\vect{\qhat}$ as the discrete solution when using isogeometric discretizations. 
The discrete director, $\vect{d}_h$, and strain/stress measures follow directly from their definitions in \eqref{eq-triad} and \eqref{eq-stress}-\eqref{eq-strain}, respectively. 
We note that the discrete director field following \eqref{eq-triad} belongs to the unit sphere $S^2$ at any point along the rod.

An alternative to isogeometric discretizations is the nodal finite elements, as employed in \cite{gebhardt_2021_beam} for the studied rod formulation. 
To fulfill the required $C^1$ continuity, the rod $\phic(s,t)$ is discretized in terms of the nodal positions and directors using cubic Hermite spline functions as follows: 
\begin{align}\label{eq-discretize_nodal}
    \phic_h(s,t) = \sum_{e=1}^{n_e} \left( H_1 \vect{x}_1^e + H_2 \vect{d}_1^e + H_3 \vect{x}_2^e + H_4 \vect{d}_2^e \right) = \mat{H}(s) \, \Bar{\vect{q}} \,,
\end{align}
where 
$H_i$, $1 \leq i \leq 4$, is the standard cubic Hermite spline function, 
$\vect{x}_j^e \in \mathbb{R}^3$ and $\vect{d}_j^e \in S^2$, $j=1,2$, is the nodal position and director at the $j$-th node of the $e$-th element, $1 \leq e \leq n_e$, respectively. Here $n_e$ denotes the number of elements and $\Bar{\vect{q}} = \Bar{\vect{q}}(t) \in (\mathbb{R}^3 \times S^2)^{n_n}$ the vector of unknown time-dependent coefficients, where $n_n = n_e+1$ is the number of nodes. 
In this work, we refer to $\Bar{\vect{q}}$ as the discrete solution when using nodal finite elements. 
We note that 
the nodal directors are here independent variables and 
the definition of the director in \eqref{eq-triad} is then valid merely within the finite elements, not at the nodes. 
Since the director field following \eqref{eq-triad} belongs to the unit sphere $S^2$, preserving this structure at the nodes requires an additional constraint of unit nodal directors:
\begin{align}\label{eq:unit_d_constraint}
    \vect{d}_j^e \cdot \vect{d}_j^e = 1\,, \quad 1 \leq e \leq n_e \,, \; j = 1,2 \,.
\end{align}
In this work, we enforce this constraint using either the Lagrange multiplier method (see e.g.\cite{Fernandez-Mendez2004,Codina2015}), as employed in \cite{gebhardt_2021_beam}, or the penalty method (see e.g. \cite{Pasch2021,Lu2019,Schillinger2016}), which we discuss more in details in the next subsection.

We note that 
with $\vect{d}_1^j \in S^2$, $1 \leq e \leq n_e$, $j = 1,2$, the resulting discrete 
solution $\Bar{\vect{q}}$ 
lives in multiple copies of the manifold $\mathbb{R}^3 \times S^2$.  
Notably, with unit nodal director, we obtain the following at any $i$-th node:
\begin{align}
    \phic_h^\prime(s_i,t) = \vect{d}_i^e\,,
\end{align}
and thus $\abss{\phic_h^\prime(s_i,t)} = 1$ 
and $\vect{d}_h((s_i,t)) = \phic_h^\prime(s_j,t)$, following \eqref{eq-triad}. 
This leads to zero axial stress at the corresponding node: 
\begin{align}
    \vect{n}_h = EA \, \vect{\epsilon}_h \,=\, EA \, \left( \phic_h^\prime \,-\, \vect{d}_h \right) = \vect{0}\,.
\end{align}
To tackle this, one can either employ isogeometric discretizations \eqref{eq-discretize_iga} or allow arbitrary length for the nodal directors, i.e. neglecting the unit nodal director constraint \eqref{eq:unit_d_constraint}. 
This leads to a discrete solution in multiple copies of the Euclidean space $\mathbb{R}^3$ that is the same  
as using isogeometric discretizations. 
In this work, we consider all three variants of the spatial discretization: 
the isogeometric discretization scheme, 
the nodal discretization scheme with and without unit nodal director constraint. 
In the next subsection, we discuss the resulting semi-discrete formulation for each of these three formulations and the enforcement of the unit nodal director constraint.

\subsection{Semi-discrete formulations}\label{sec:semi_discrete_form}

Using isogeometric discretizations, 
we recall from \cite{nguyen_rod_2024} the semi-discrete formulation, obtained after introducing \eqref{eq-discretize_iga} into the weak form \eqref{w-eom}:
\begin{align}\label{semi-deom_iga}
    & \text{Find } \vect{\qhat}(t) \,\in\, \mathbb{R}^{3\,\m}, \, t\,\in\,(0,\,T], \text{ such that}: \nonumber \\
    & \int_0^S \, \delta \vect{\qhat} \, \cdot \, \left(\, \mat{M} (\vect{\qhat}) \, \nabla_{\dot{\vect{\qhat}}} \, \dot{\vect{\qhat}} \,+\, \mat{B} (\vect{\qhat})^T \, \vect{\sigma}_h -\, \mat{\bspline}^T \, \vect{f}^{\text{ext}} \,\right) \, \mathrm{d} \, s \,=\, \mat{0} \quad \forall \, \delta \vect{\qhat} \,\in\, \mathbb{R}^{3\,\m} \,.
\end{align}
Here, the mass matrix $\mat{M}$ and the matrix $\mat{B}$ result from the operators \eqref{mass-op} and \eqref{B-op}, respectively, 
$\vect{\sigma}_h$ denotes the discrete stress measures, $\vect{\sigma}_h = \left[\vect{n}_h, \;\; \vect{m}_h\right]^T$, 
and $\vect{f}^{\text{ext}}$ the external force vector. 
For detailed derivation and explicit expressions of these matrices and vectors, we refer to our previous work \cite{nguyen_rod_2024}. 
We note that the initial conditions at $t=0$ follow from the initial conditions $\phic \,=\, \phic_0$ and $\dot{\phic} \,=\, \mat{v}_0$ discussed in the previous section.

Analogously, using nodal finite elements without considering the unit nodal director constraint leads to a similar semi-discrete formulation:
\begin{align}\label{semi-deom_nodalR3}
    & \text{Find } \Bar{\vect{\qhat}}(t) \,\in\, \mathbb{R}^{6\,n_n}, \, t\,\in\,(0,\,T], \text{ such that}: \nonumber \\
    & \int_0^S \, \delta \Bar{\vect{\qhat}} \, \cdot \, \left(\, \mat{M} (\Bar{\vect{\qhat}}) \, \nabla_{\dot{\Bar{\vect{\qhat}}}} \, \dot{\Bar{\vect{\qhat}}} \,+\, \mat{B} (\Bar{\vect{\qhat}})^T \, \boldsymbol{\sigma}_h -\, \mat{H}^T \, \vect{f}^{\text{ext}} \,\right) \, \mathrm{d} \, s 
    \,=\, \mat{0} \quad \forall \, \delta \Bar{\vect{\qhat}} \,\in\, \mathbb{R}^{6\,n_n} \,.
\end{align}
Here, the bar over the coefficient vector $\Bar{\vect{\qhat}}$ is to distinguish this from the one associated with the isogeometric scheme above.

Using nodal finite elements with consideration of the unit nodal director constraint, one can enforce this constraint \eqref{eq:unit_d_constraint} using the Lagrange multiplier method (see e.g.\cite{Fernandez-Mendez2004,Codina2015}).   
The resulting semi-discrete formulation is then:
\begin{align}\label{semi-deom_nodalR3S2-strong}
    & \text{Find } \Bar{\vect{\qhat}}(t) \,\in\, (\mathbb{R}^3 \times S^2)^{n_n}, \, t\,\in\,(0,\,T], \text{ such that}: \nonumber \\
    & \int_0^S \, \delta \Bar{\vect{\qhat}} \, \cdot \, \left(\, \mat{M} (\Bar{\vect{\qhat}}) \, \nabla_{\dot{\Bar{\vect{\qhat}}}} \, \dot{\Bar{\vect{\qhat}}} \,+\, \mat{B} (\Bar{\vect{\qhat}})^T \, \boldsymbol{\sigma}_h -\, \mat{H}^T \, \vect{f}^{\text{ext}} \,\right) \, \mathrm{d} \, s \nonumber \\
    & \qquad \qquad \qquad +\, \delta \Bar{\vect{\qhat}} \cdot \mat{J}^T (\Bar{\vect{\qhat}}) \,\vect{\lambda} \,+\, \delta \vect{\lambda} \cdot \vect{\Psi}
    \,=\, \mat{0} \quad 
    \forall \, \delta \Bar{\vect{\qhat}} \,\in\, (\mathbb{R}^3 \times S^2)^{n_n} \,, \;
    \delta \vect{\lambda} \,\in\, \mathbb{R}^{n_n} \,,
\end{align}
where
$\vect{\lambda}$ and $\vect{\Psi}$ are the vector of unknown Lagrange multipliers and the vector of unit director constraint at each node, that are:
\begin{align}
    & \vect{\lambda} = \left[\lambda_1 \; \ldots \; \lambda_{n_n} \right]^T\,, \\
    & \vect{\Psi} = \left[\vect{d}_1 \cdot \vect{d}_1 - 1 \; \ldots \; \vect{d}_{n_n} \cdot \vect{d}_{n_n} - 1 \right]^T\,,
\end{align}
respectively. The matrix $\mat{J}= \mat{J}(\Bar{\vect{\qhat}})$ results from the variation of $\vect{\Psi}$ that is:
\begin{align}\label{eq:variation_of_constraint}
    & \delta \boldsymbol{\Psi} = 2\,\left[\delta \vect{d}_1 \cdot \vect{d}_1 \; \ldots \; \delta \vect{d}_{n_n} \cdot \vect{d}_{n_n} \right]^T = \mat{J} \, \delta \, \Bar{\vect{\qhat}}\,.
\end{align}
For the explicit expression of $\mat{J}$, we refer to Appendix \ref{sec:linearized_constraint_vector}. 
We note that this approach using the Lagrange multiplier method requires an additional variable field of the Lagrange multipliers in $\vect{\lambda}$. 
In this work, to eliminate this variable field 
we employ the nullspace matrix $\mat{D} = \mat{D}(\Bar{\vect{\qhat}})$ of the matrix $\mat{J}$, i.e. $\mat{J} \, \mat{D} = \mat{0}$. 
In particular, we 
consider $\delta \Bar{\vect{\qhat}} = \left(\mat{D} \delta \tilde{\vect{\qhat}}\right)$ and obtain the following semi-discrete formulation: 
\begin{align}\label{semi-deom_nodalR3S2-strong-reduced}
    & \delta \tilde{\vect{\qhat}} \, \cdot \, \left[ \mat{D}(\Bar{\vect{\qhat}})^T \; \int_0^S \, \left(\, \mat{M} (\Bar{\vect{\qhat}}) \, \nabla_{\dot{\Bar{\vect{\qhat}}}} \, \dot{\Bar{\vect{\qhat}}} \,+\, \mat{B} (\Bar{\vect{\qhat}})^T \, \boldsymbol{\sigma}_h -\, \mat{H}^T \, \vect{f}^{\text{ext}} \,\right) \, \mathrm{d} \, s \right] + 
    \delta \vect{\lambda} \cdot \vect{\Psi} \,=\, \mat{0} \,.
\end{align}
For the explicit expression of $\mat{D}$, we refer to Appendix \ref{sec:linearized_constraint_vector}.

Alternatively, one can enforce the unit nodal director constraint \eqref{eq:unit_d_constraint} using the penalty method (see e.g. \cite{Pasch2021,Lu2019,Schillinger2016}). 
Inspired by the work of\cite{Pasch2021}, we remove the problem-denpendency in choosing the penalty factor  
and obtain the unit consistency of the penalty term via scaling the penalty factor by $2EI/L$. 
The corresponding semi-discrete formulation is then:
\begin{align}\label{semi-deom_nodalR3S2-weak}
    & \text{Find } \Bar{\vect{\qhat}}(t) \,\in\, \mathbb{R}^{6\,n_n}, \, t\,\in\,(0,\,T], \text{ such that}: \nonumber \\
    & \int_0^S \, \delta \Bar{\vect{\qhat}} \, \cdot \, \left(\, \mat{M} (\Bar{\vect{\qhat}}) \, \nabla_{\dot{\Bar{\vect{\qhat}}}} \, \dot{\Bar{\vect{\qhat}}} \,+\, \mat{B} (\Bar{\vect{\qhat}})^T \, \boldsymbol{\sigma}_h -\, \mat{H}^T \, \vect{f}^{\text{ext}} \,\right) \, \mathrm{d} \, s \, 
    +\, \beta \, \frac{2EI}{L} \delta \Bar{\vect{\qhat}} \cdot \mat{J}^T (\Bar{\vect{\qhat}}) \,\vect{\Psi}
    \,=\, \mat{0} \quad 
    \forall \, \delta \Bar{\vect{\qhat}} \,\in\, \mathbb{R}^{6\,n_n} \,,
\end{align}
where $\beta$ is the penalty factor. 
Here, we choose a scaling factor of $2EI/L$ which is analytically determined for Timoshenko beams when enforcing Dirichlet boundary conditions using the penalty method in \cite{Pasch2021}. 
This choice is based on the fact that the unit nodal director constraint involves the nodal director degrees of freedom, which are associated with the rotation of the rod axis. 
We note that this approach using the penalty method does not strictly lead to unit nodal director since an infinite penalty factor cannot be employed for numerical computations. 
Hence, the resulting discrete solution $\Bar{\vect{\qhat}}$ is strictly in $(\mathbb{R}^3)^{2n_n}$ instead of the manifold $(\mathbb{R}^3 \times S^2)^{n_n}$. 

\begin{table}[ht]
    \centering
    \begin{tabularx}{1\linewidth}{|X | X | >{\centering\arraybackslash}X| >{\centering\arraybackslash}X|}
        \toprule
        \multicolumn{2}{|l|}{\textbf{Discretization scheme}} & \textbf{Semi-discrete formulation} & \textbf{Discrete solution} \\ 
        \hline
        \multicolumn{2}{|l|}{Isogeometric discretizations} & Equation \eqref{semi-deom_iga} & $\mathbb{R}^{3\,\m}$ \\
        \cline{3-4}
        \multicolumn{4}{|l|}{\textit{Note that the axial stress resultant is not constrained to zero at any point and the discrete director field lives in $S^2$ at any}} \\
        \multicolumn{4}{|l|}{\textit{point of the discrete rod configuration.}} \\
        \hline
        \multicolumn{2}{|l|}{Nodal discretization scheme without unit nodal director constraint} & Equation \eqref{semi-deom_nodalR3} & $\mathbb{R}^{6\,n_n}$ \\
        \cline{3-4}
        \multicolumn{4}{|l|}{\textit{Note that the nodal axial stress is not constrained to zero, however, nodal directors and director defined within elements}} \\
        \multicolumn{4}{|l|}{\textit{live in different spaces: $\mathbb{R}^3$ and $S^2$, respectively.}} \\
        \hline 
        \multicolumn{1}{|l|}{Nodal discretization} & \multicolumn{1}{l|}{enforcement using Lagrange multiplier} & Equation \eqref{semi-deom_nodalR3S2-strong} & $(\mathbb{R}^3 \times S^2)^{n_n}$ \\
        \multicolumn{1}{|l|}{scheme with unit nodal} & \multicolumn{1}{l|}{method} & & \\
        \cline{2-4}
        \multicolumn{1}{|l|}{director constraint} & \multicolumn{1}{l|}{enforcement with reduced equations us-} & Equation \eqref{semi-deom_nodalR3S2-strong-reduced} & $(\mathbb{R}^3 \times S^2)^{n_n}$ \\
        & \multicolumn{1}{l|}{ing Lagrange multiplier and nullspace} & & \\
        & \multicolumn{1}{l|}{methods} & & \\
        \cline{2-4}
        & \multicolumn{1}{l|}{enforcement using penalty method} & Equation \eqref{semi-deom_nodalR3S2-weak} & (strictly)$^1$ $\mathbb{R}^{6\,n_n}$ \\ 
        \cline{2-4}
        \multicolumn{4}{|l|}{\textit{Note that nodal axial stress is zero, however, the discrete director field lives in $S^2$ at any point of the discrete rod configura-}} \\
        \multicolumn{4}{|l|}{\textit{tion.}} \\
        \hline
        \multicolumn{4}{|l|}{\textit{$^1$ With a penalty factor $\beta \rightarrow \infty$, the discrete solution space becomes $(\mathbb{R}^3 \times S^2)^{n_n}$.}} \\
        \bottomrule
    \end{tabularx}
    \caption{Different semi-discrete formulations and discrete solution spaces for the shear- and torsion-free nonlinear rod \cite{gebhardt_2021_beam}, using isogeometric and nodal discretizations.}\label{tab:semi-discrete-forms}
\end{table}

\begin{remark}
    We note that the factor $\beta$ in Equation \eqref{semi-deom_nodalR3S2-weak} is not the intensity factor\footnote{An intensity factor of $\beta$ leads to an error of $(100/\beta)$\% between the solution obtained with the penalty and the Lagrange multiplier method \cite{Pasch2021}.} as defined in \cite{Pasch2021} since we neither consider the Timoshenko beam formulation nor apply the penalty method to enforce the Dirichlet boundary conditions. 
    Nevertheless, in this work, $\beta$ is a scaled intensity factor and provides an estimate of the solution error. One can choose a sufficiently large value of $\beta$ to achieve a small solution error. 
    An analytical determination of the scaling factor, as developed in \cite{Pasch2021}, is out of scope of this work and is therefore left for future investigation.
\end{remark}

In this work, we consider and compare all five semi-discrete formulations discussed above, for which we discuss the resulting matrix equations in the next section. 
We give an overview of these five formulations in Table \ref{tab:semi-discrete-forms}.

\section{Matrix equations}\label{sec:matrix_form}
  
In this section, we review existing and derive matrix equations of five semi-discrete formulations studied in this work, as listed in Table \ref{tab:semi-discrete-forms}. 
One obtains the corresponding matrix equations after employing the implicit time integration scheme reviewed in Section \ref{sec:pre_time_integration} and linearizing the resulting nonlinear residual. 
We start with reviewing the equations obtained with isogeometric discretizations.

\begin{remark}
    In this work, we focus on dynamics cases governed by the semi-discrete equations discussed in the previous section. For static cases, we refer to \cite{nguyen_rod_2024} for more details on the corresponding semi-discrete and matrix equations.
\end{remark}

\subsection{Isogeometric discretizations}

Using the implicit time integration scheme \cite{gebhardt_2021_beam,nguyen_rod_2024} (see also Section \ref{sec:pre_time_integration}) requires to evaluate the semi-discrete formulation \eqref{semi-deom_iga} at the time instance $t_{n+\frac{1}{2}}$. 
This can be then approximated using the midpoint and trapezoidal rules, resulting in a system of nonlinear equations depending on the known discrete solution at $t_{n}$ and unknown solution at $t_{n+1}$. 
Using Newton-Raphson method, we solve the resulting equations for the discrete solution at $t_{n+1}$. 
When using isogeometric discretizations, we refer to \cite{nguyen_rod_2024} for details of the approximation for each term of the semi-discrete formulation \eqref{semi-deom_iga}. 
We recall here the resulting linearized matrix equations:
\begin{align}\label{eq:matrix_eq_iga}
    \mat{A}\left(\vect{\qhat}_{n+1}^{k-1}, \vect{\qhat}_{n}, \dot{\vect{\qhat}}_{n}\right) \, \Delta \vect{\qhat}_{n+1}^k = \vect{F}^{\text{ext}}_{n+\frac{1}{2}} - \vect{F}\left(\vect{\qhat}_{n+1}^{k-1}, \vect{\qhat}_{n}, \dot{\vect{\qhat}}_{n} \right) \,,
\end{align}
where $k$ denotes the iteration step of the Newton-Raphson scheme 
and the subscript refers to the time instance at which the term is evaluated. 
For the derivation and computation of the system matrix $\mat{A}$ and the force vector $\vect{F}$, we refer to \cite{nguyen_rod_2024}.

We note that 
the matrix $\mat{A}$ is a sparse and symmetric matrix. 
As discussed in \cite{nguyen_rod_2024}, the number of degrees of freedom (dofs) can be estimated based on the number of elements and the chosen spline basis functions. 
A discretization with $n_e$ elements using $C^r$ B-splines of degree $p$ leads to $3 [ne (p-r) + r + 1]$ dofs, assuming that the splines are defined on an open knot vector with interior knots repeated $(p-r)$-times. 
Also reported in \cite{nguyen_rod_2024}, using isogeometric discretizations can lead to responses with high-frequency contents in some cases, which can be reduced using the strong approach of outlier removal \cite{hiemstra_outlier_2021}. 
This outlier removal approach effectively and entirely removes the spurious outlier modes that correspond to the highest frequencies in the discrete spectrum, without compromising the accuracy. 
It requires the computation of the extraction operator $\mat{C}$ that builds constraints for outlier removal into the space of spline basis functions \cite{nguyen_rod_2024,hiemstra_outlier_2021}. 
$\mat{C}$ only depends on the spline space and required constraints, and hence is constant and can be computed only once before time integration. 
In each iteration at each time step, the system matrix is globally multiplied from left and right by $\mat{C}$, as well as the right-hand side from the left. \eqref{eq:matrix_eq_iga} becomes:
\begin{align}\label{eq:matrix_eq_iga_outlier_removal}
     \mat{C}^T \mat{A} \mat{C} \, \Delta \vect{\qhat}_{n+1}^k = \mat{C}^T \left( \vect{F}^{\text{ext}}_{n+\frac{1}{2}} - \vect{F}\right) \,.
\end{align}
The resulting system matrix on the left-hand side remains sparse, symmetric, and has smaller dimensions than the original one.

\subsection{Nodal discretization scheme without unit nodal director constraint}

Neglecting the unit nodal director constraint when using the nodal discretization scheme, the resulting semi-discrete formulation \eqref{semi-deom_nodalR3} is similar to \eqref{semi-deom_iga} obtained with isogeometric discretizations. 
Analogously, employing the same implicit time integration scheme \cite{gebhardt_2021_beam,nguyen_rod_2024} (see also Section \ref{sec:pre_time_integration}) leads to the following matrix equations:
\begin{align}\label{eq:matrix_eq_nodalR3}
    \Bar{\mat{A}}\left(\Bar{\vect{\qhat}}_{n+1}^{k-1}, \Bar{\vect{\qhat}}_{n}, \dot{\Bar{\vect{\qhat}}}_{n}\right) \, \Delta \Bar{\vect{\qhat}}_{n+1}^k = \Bar{\vect{F}}^{\text{ext}}_{n+\frac{1}{2}} - \Bar{\vect{F}}\left(\Bar{\vect{\qhat}}_{n+1}^{k-1}, \Bar{\vect{\qhat}}_{n}, \dot{\Bar{\vect{\qhat}}}_{n} \right) \,,
\end{align}
where $k$ denotes the iteration step of the Newton-Raphson scheme. For the derivation and computation of the system matrix $\Bar{\mat{A}}$ and the force vector $\Bar{\vect{F}}$, we refer to \cite{gebhardt_2021_beam}. 
Similarly to \eqref{eq:matrix_eq_iga}, $\Bar{\mat{A}}$ is a sparse and symmetric matrix.

\begin{figure}[ht]
	\centering
	\subfloat[Cubic $C^1$ B-splines]{\includegraphics[width=0.5\textwidth]{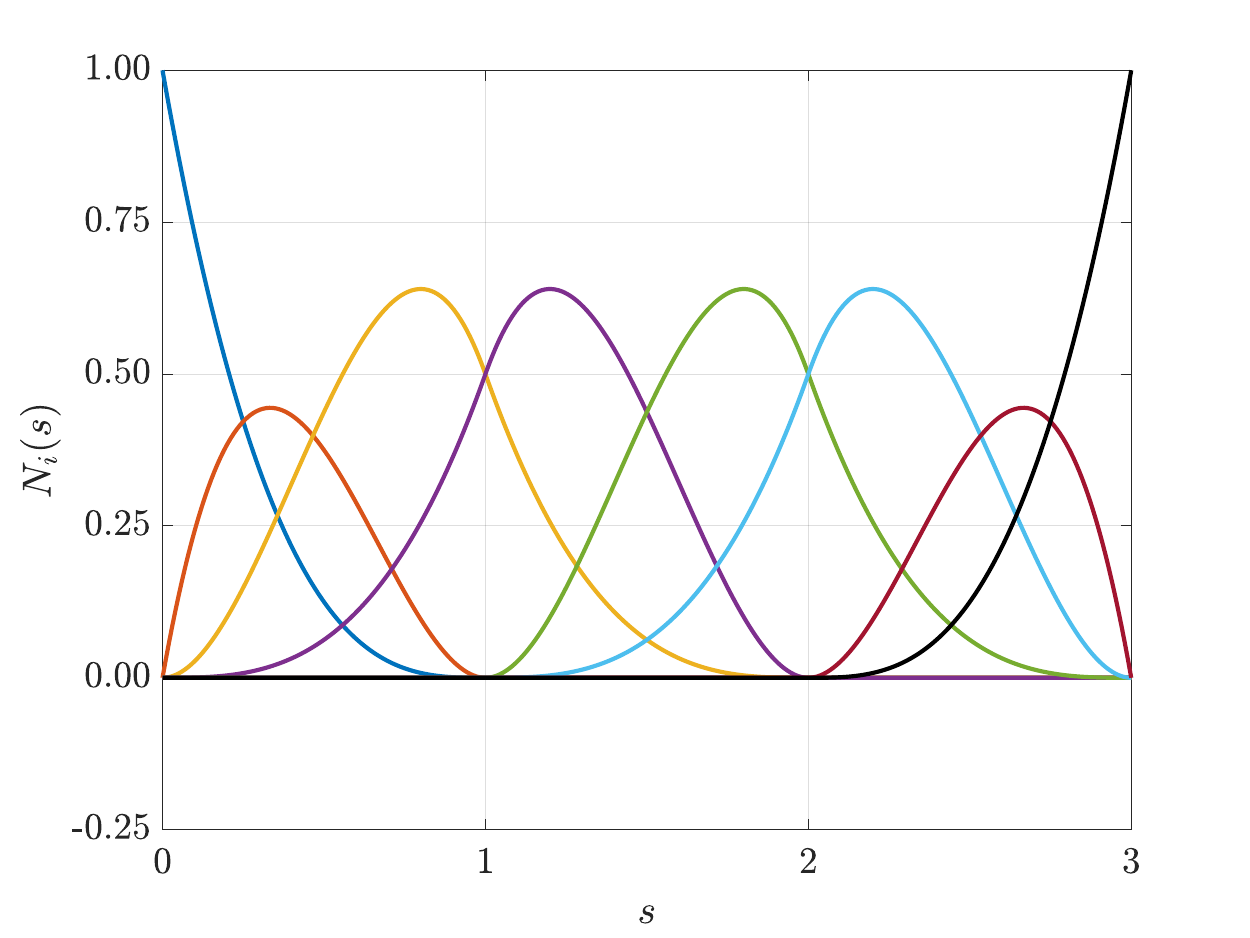}} 
    \subfloat[Cubic Hermite splines]{\includegraphics[width=0.5\textwidth]{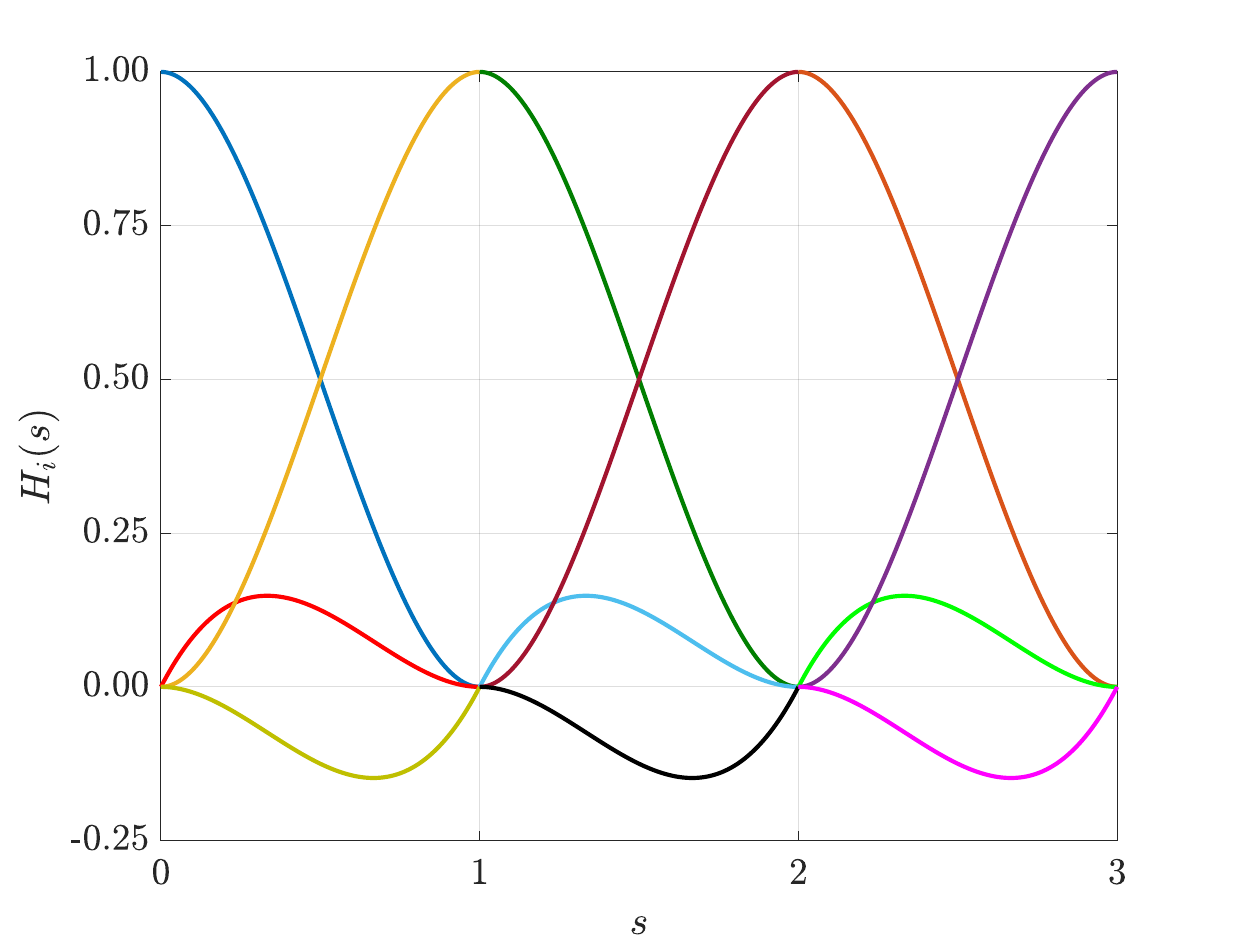}}
    \vspace{0.3cm}
    
	\caption{Cubic $C^1$ B-splines and cubic Hermite splines defined on an interval uniformly discretized in three elements. These two bases span the same function space.}\label{fig:spline_hermite}
\end{figure}

When using nodal discretization scheme, the basis functions are the standard cubic Hermite splines. 
A discretization with $n_e$ elements leads to $6(n_e+1)$ dofs, i.e. the size of the vector $\Delta \Bar{\vect{\qhat}}$ consisting of unknown coefficients. 
We note that 
the cubic Hermite spline and the cubic $C^1$ B-splines span the same function space (see also Fig. \ref{fig:spline_hermite}).

\subsection{Nodal discretization scheme with unit nodal director constraint}

We now discuss the resulting matrix equations when enforcing the unit nodal director constraint \eqref{eq:unit_d_constraint} using different approaches. 
We focus here on the resulting equations for comparison purposes  
and refer 
to Appendix \ref{sec:linearized_constraint_vector} 
for 
the derivation of associated matrices and further technical details.

    \subsubsection{Enforced constraint with Lagrange multipliers}

    Consider the semi-discrete formulation \eqref{semi-deom_nodalR3S2-strong}. 
    Employing the implicit time integration scheme reviewed in Section \ref{sec:pre_time_integration} and linearizing the resulting nonlinear residual leads to the following matrix equations at the $k$-th iteration:
    \begin{align}\label{eq:matrix_eq_nodalR3S2-strong}
        \underbrace{\begin{bmatrix}
            \Bar{\mat{A}} + \mat{A}_c & \mat{J}^T_{n+\frac{1}{2}} \\
            \mat{J}_{n+1} & \mat{0}
        \end{bmatrix}}_{\Tilde{\mat{A}}} \, \begin{bmatrix}
            \Delta \bar{\vect{\qhat}}_{n+1}^k \\ \Delta \vect{\lambda}_{n+\frac{1}{2}}^k
        \end{bmatrix} = \begin{bmatrix}
            \Bar{\vect{F}}^{\text{ext}}_{n+\frac{1}{2}} - \Bar{\vect{F}} - \mat{J}^T_{n+\frac{1}{2}} \vect{\lambda}_{n+\frac{1}{2}} \\
            - \vect{\Psi}_{n+1}
        \end{bmatrix}\,,
    \end{align}
    where the matrix 
    $\Bar{\mat{A}} = \Bar{\mat{A}}\left(\bar{\vect{\qhat}}_{n+1}^{k-1}, \bar{\vect{\qhat}}_{n}, \dot{\bar{\vect{\qhat}}}_{n} \right)$ is the same as that in \eqref{eq:matrix_eq_nodalR3}, 
    $\mat{A}_c = \mat{A}_c\left(\vect{\lambda}_{n+\frac{1}{2}}^{k-1} \right)$ is the contribution of the unit nodal director constraint to the system matrix, i.e. the linearization of the term $\delta \Bar{\vect{\qhat}} \cdot \mat{J}^T (\Bar{\vect{\qhat}}) \,\vect{\lambda}$ evaluated at $t_{n+\frac{1}{2}}$. 
    For the derivation of $\mat{A}_c$, we refer to Appendix \ref{sec:linearized_constraint_vector}.

    We note that the two blocks of $\mat{J}^T$ and $\mat{J}$ of $\Tilde{\mat{A}}$ are evaluated at two different time instances due to the holomonic type of the unit nodal director constraint, as discussed in \cite[Remark~2, p.~3834]{gebhardt_2021_beam}. 
    Hence, $\Tilde{\mat{A}}$ on the left-hand side of \eqref{eq:matrix_eq_nodalR3S2-strong} is not a symmetric matrix, however, is a sparse matrix. 
    \eqref{eq:matrix_eq_nodalR3S2-strong} has the form of a saddle-point problem (SPP), for which we briefly recall the equivalence of the necessary and sufficient conditions for unique solution \cite[p.~142]{brezzi_mixedFEM2013}:
    \begin{itemize}
        \item The matrix $\left(\Bar{\mat{A}} + \mat{A}_c \right)$ is symmetric positive semi-definite, and
        \item The matrix $\mat{J}$ is full rank.
    \end{itemize}
    In general, these conditions are fulfilled for the studied rod formulation due to the definition of these matrices (see \cite{nguyen_rod_2024} for $\Bar{\mat{K}}$ and Appendix \ref{sec:linearized_constraint_vector} for $\mat{A}_c$ and $\mat{J}$). 
    We note that the number of degrees of freedom (dofs) is at most $(n_e+1)$ dofs more than that of     
    \eqref{eq:matrix_eq_nodalR3}, i.e. $7(n_e+1)$ dofs, due to at most $(n_e+1)$ additional unknown Lagrange multipliers. 
    In cases of a clamped boundary condition or prescribed director, it requires a smaller number of unknown Lagrange multipliers. 
    We note that here, 
    we refer to the size of the vector     
    $\left[\Delta \bar{\vect{\qhat}}_{n+1}^k \; \Delta \vect{\lambda}_{n+\frac{1}{2}}^k \right]^T$ as the number of dofs, i.e. the number of unknown variables, for which we solve at each iteration.

    \subsubsection{Enforced constraint with Lagrange multiplier and nullspace methods}

    In this work, to eliminate the additional variable field of Lagrange multipliers in \eqref{semi-deom_nodalR3S2-strong}, i.e. to reduce the dimension of \eqref{eq:matrix_eq_nodalR3S2-strong} to the same as that of \eqref{eq:matrix_eq_nodalR3}, 
    we employ the nullspace matrix $\mat{D}$ of the matrix $\mat{J}$, as discussed in the previous section. 
    Consider the semi-discrete formulation \eqref{semi-deom_nodalR3S2-strong-reduced}, employing the implicit time integration scheme reviewed in Section \ref{sec:pre_time_integration} and linearizing the resulting nonlinear residual leads to the following matrix equations at the $k$-th iteration:
    \begin{align}\label{eq:matrix_eq_nodalR3S2-strong-reduced}
        \underbrace{\begin{bmatrix}
            \mat{D}^T_{n+\frac{1}{2}} \Bar{\mat{A}} + \mat{A}_D \\
            \mat{J}_{n+1}
        \end{bmatrix}}_{\hat{\mat{A}}} \, \Delta \bar{\vect{\qhat}}_{n+1}^k = \begin{bmatrix}
            \mat{D}^T_{n+\frac{1}{2}} \left(\Bar{\vect{F}}^{\text{ext}}_{n+\frac{1}{2}} - \Bar{\vect{F}}\right) \\
            - \vect{\Psi}_{n+1}
        \end{bmatrix}\,,
    \end{align}
    where the counterpart $\mat{K}_D$ is the contribution of the nullspace matrix to the system matrix, i.e. the linearization of $\mat{D}^T_{n+\frac{1}{2}}$. 
    For its derivation, we refer to Appendix \ref{sec:linearized_constraint_vector}. 
    We note that $\mat{D}^T_{n+\frac{1}{2}}$ depends on the nodal director from the previous iteration (see also Appendix \ref{sec:linearized_constraint_vector}) and hence needs to be reassembled in each iteration and time step. 
    Furthermore,     
    the matrix multiplication by $\mat{D}^T_{n+\frac{1}{2}}$ is performed globally, also at each time step and iteration. 
    The resulting matrix $\hat{\mat{K}}$ on the left-hand side of \eqref{eq:matrix_eq_nodalR3S2-strong-reduced} is a sparse but not a symmetric matrix. 
    The number of degrees of freedom is now the same as that of \eqref{eq:matrix_eq_nodalR3}.

    \subsubsection{Enforced constrained with penalty method}

    An alternative approach is to enforce the unit nodal director constraint \eqref{eq:unit_d_constraint} using the penalty method. 
    Consider the semi-discrete formulation \eqref{semi-deom_nodalR3S2-weak}, 
    employing the implicit time integration scheme reviewed in Section \ref{sec:pre_time_integration} and linearizing the resulting nonlinear residual leads to the following matrix equations at the $k$-th iteration:
    \begin{align}\label{eq:matrix_eq_nodalR3S2-weak}
        \underbrace{\left( \Bar{\mat{A}} + \mat{A}_{\beta,n+1} \right)}_{\check{\mat{A}}} 
        \, \Delta \Bar{\vect{\qhat}}_{n+1}^k = \Bar{\vect{F}}^{\text{ext}}_{n+\frac{1}{2}} - \Bar{\vect{F}} - \beta \, \frac{2EI}{L} \, \mat{J}^T_{n+1} \vect{\Psi}_{n+1} \,,
    \end{align}
    where $\mat{A}_{\beta,n+1}$ is the contribution of the penalty term to the system matrix, i.e. the linearization of the penalty term $\left( \beta \frac{2EI}{L} \, \delta \Bar{\vect{\qhat}} \cdot \mat{J}^T (\Bar{\vect{\qhat}}) \,\vect{\Psi} \right)$ evaluated at the time instance $t_{n+1}$. 
    For its derivation, we refer to Appendix \ref{sec:linearized_constraint_vector}. 
    We note that the resulting matrix $\check{\mat{A}}$ on the left-hand side of \eqref{eq:matrix_eq_nodalR3S2-weak} is a sparse and symmetric matrix. 
    The number of degrees of freedom is the same as that of \eqref{eq:matrix_eq_nodalR3}.

\section{Computational cost}\label{sec:computational_cost}
  
In this section, we discuss the computation cost corresponding to the five formulations discussed in the previous two sections, listed in Table \ref{tab:semi-discrete-forms}. 
Particularly, we compare the sparsity, bandwidth, and symmetry of the system matrix. 
We also discuss the number of degrees of freedom (dofs) of each formulation. 
In Table \ref{tab:matrix-forms}, we give an overview of these properties for the five studied formulations. 
For the notation simplicity, we refer to the studied formulations using the abbreviations given in italics in Table \ref{tab:matrix-forms}.

\begin{table}[ht]
    \centering
    \begin{tabularx}{1\linewidth}{|X | X | >{\centering\arraybackslash}X| >{\centering\arraybackslash}X| >{\centering\arraybackslash}X|}
        \toprule
        \multicolumn{2}{|l|}{\textbf{Discretization scheme}} & \textbf{Matrix equations} & \textbf{System matrix} & \textbf{Number of dofs}$^1$ \\ 
        \hline
        \multicolumn{2}{|l|}{Isogeometric discretizations (\textit{\iga})} & Equation$^2$ \eqref{eq:matrix_eq_iga}, \eqref{eq:matrix_eq_iga_outlier_removal}$^\star$ & sparse, symmetric & $3 [ne (p-r) + r + 1]$ \\
        \hline
        \multicolumn{2}{|l|}{Nodal discretization scheme without unit nodal director} & Equation \eqref{eq:matrix_eq_nodalR3} & sparse, symmetric & $6(n_e+1)$ \\
        \multicolumn{2}{|l|}{constraint (\textit{\nodal})} &  &  & \\
        \hline 
        \multicolumn{1}{|l|}{Nodal discretization} & \multicolumn{1}{l|}{enforcement using Lagrange mul-} & Equation \eqref{eq:matrix_eq_nodalR3S2-strong}& sparse, non-& $7(n_e+1)$ \\
        \multicolumn{1}{|l|}{scheme with unit} & \multicolumn{1}{l|}{tiplier method (\textit{\nodalSaddle})} &  & symmetric & \\
        \cline{2-5}
        \multicolumn{1}{|l|}{nodal director con-} & \multicolumn{1}{l|}{enforcement with reduced equa-} & Equation \eqref{eq:matrix_eq_nodalR3S2-strong-reduced}$^{\star\star}$ & sparse, non- & $6(n_e+1)$ \\
        \multicolumn{1}{|l|}{straint} & \multicolumn{1}{l|}{tions using Lagrange multiplier} & & symmetric & \\
        & \multicolumn{1}{l|}{and nullspace methods} & & & \\
        & \multicolumn{1}{l|}{(\textit{\nodalSaddleRed})} & & & \\
        \cline{2-5}
        & \multicolumn{1}{l|}{enforcement using penalty method} & Equation \eqref{eq:matrix_eq_nodalR3S2-weak} & sparse, symmetric & $6(n_e+1)$ \\
        & \multicolumn{1}{l|}{(\textit{\nodalPenalty})} &  & & \\   
        \hline
        \multicolumn{5}{|l|}{$^\star$: Global matrix multiplication is required. The multiplier is a constant matrix.} \\
        \multicolumn{5}{|l|}{$^{\star\star}$: Global matrix multiplication is required. The multiplier is reassembled at each iteration and time step.} \\
        \multicolumn{5}{|l|}{\textit{$^1$ i.e. the size of the unknown increment vector on the left-hand side, for which we solve at each iteration.}} \\
        \multicolumn{5}{|l|}{\textit{$^2$ Using the strong approach of outlier removal \cite{hiemstra_outlier_2021} leads to a smaller number of dofs than $3 [ne (p-r) + r + 1]$, depending on}} \\
        \multicolumn{5}{|l|}{\textit{the type of Dirichlet boundary conditions.}} \\
        \bottomrule
    \end{tabularx}
    \caption{Different matrix equations for the studied semi-discrete formulations listed in Table \ref{tab:semi-discrete-forms}, using either isogeometric or nodal discretizations.}\label{tab:matrix-forms}
\end{table}

Focusing on the number of degrees of freedom (dofs) when using the nodal discretization scheme, we see that enforcing the unit nodal director constraint using the Lagrange multiplier method (\textit{\nodalSaddle}) leads to the highest number of dofs due to the additional unknown Lagrange multipliers. 
This can be reduced to the same number of dofs $6(n_e+1)$ by eliminating the Lagrange multipliers using the nullspace method (\textit{\nodalSaddleRed}), or enforcing the same constraint using the penalty method (\textit{\nodalPenalty}), or neglecting the unit nodal director constraint (\textit{\nodal}) (see Table \ref{tab:matrix-forms}). 
Focusing on the number of degrees of freedom (dofs) when using isogeometric discretizations, employing quadratic $C^1$ splines leads to $3(n_e+2)$ dofs, less than $6(n_e+1)$ when using nodal discretizations with the same number of elements. 
Using cubic $C^1$ splines that span the same function space as the cubic Hermite splines employed for the nodal scheme, leads to $6(n_e+1)$ dofs, the same as 
the smallest number of dofs when using the nodal scheme. 
Using cubic $C^2$ splines leads to either the same dofs of $6(n_e+1)$ when using one element, i.e. $n_e=1$, or less than $6(n_e+1)$ with $n_e>1$. 
Using splines of higher polynomial degrees $p$ and higher order of continuity $r$ may lead to more than $6(n_e+1)$ dofs, except the cases with a very large number of elements. 
In summary, 
the isogeometric discretization scheme enables quadratic basis functions that are one order lower than Hermite splines and 
a smaller number of dofs than using the nodal scheme, particularly when using a significantly large number of elements. 
We note that for the studied rod formulation, a significantly large number of elements is generally not required and thus the number of dofs is not decisive for the difference in the computational cost when using isogeometric or nodal discretizations.

\begin{figure}
	\centering
    \def\svgwidth{0.83\textwidth}
    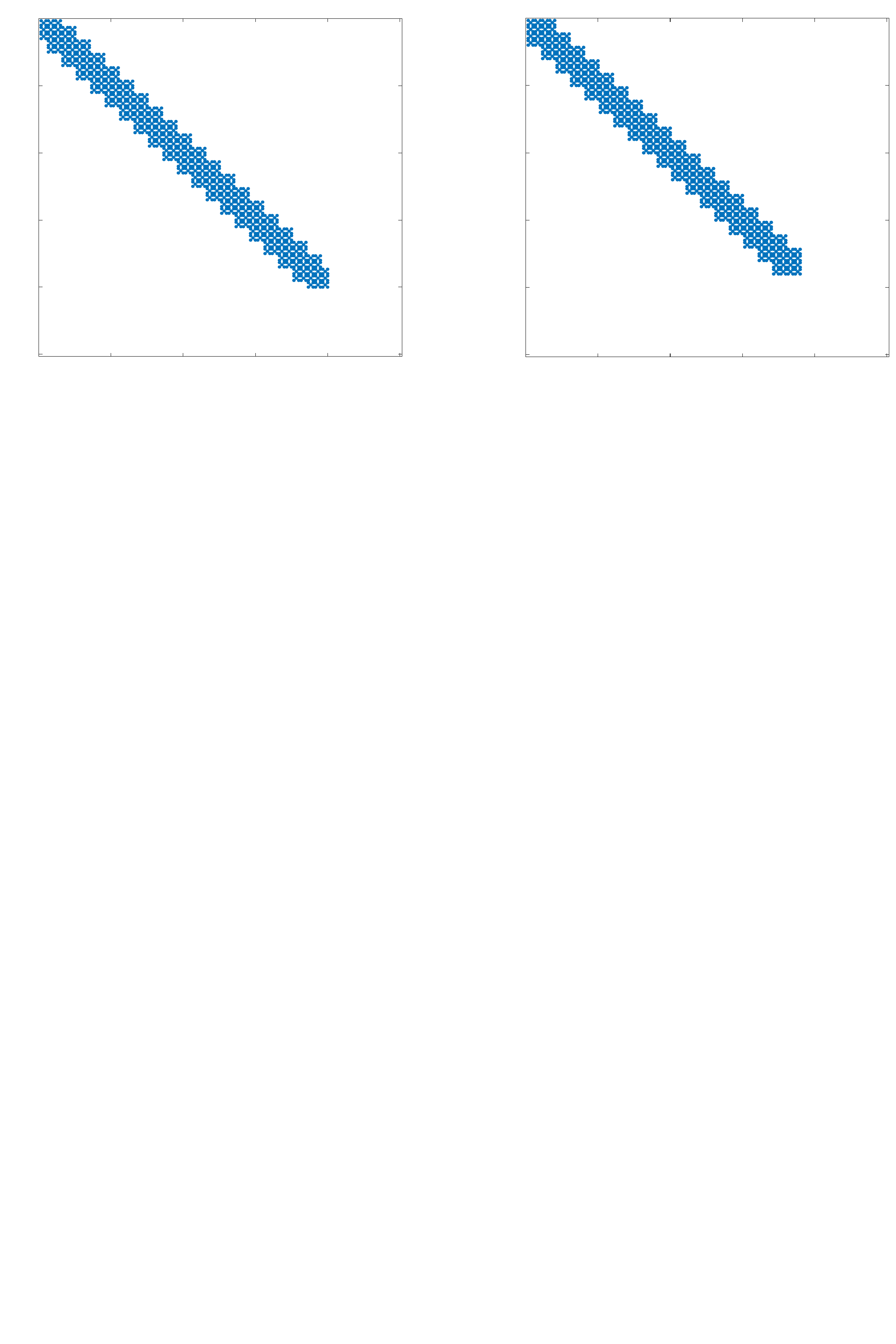
    
	\caption{Band structure of the system matrix using isogeometric and nodal discretization schemes with different approaches.} \label{fig:matrix_structure}
\end{figure}

Focusing on the sparsity of the system matrix, we see that we obtain a sparse matrix in all cases. 
Moreover, as discussed above, a significantly large number of elements and thus also the number of dofs is generally not required, using either isogeometric or nodal discretizations requires similar memory storage for the system matrix. 
Another essential factor that affects the computational cost per iteration is the symmetry of the system matrix since solving a symmetric system generally requires less effort than solving a non-symmetric one when using a standard solver. 
As discussed in the previous section and summarized in Table \ref{tab:matrix-forms}, 
we observe that 
using the nodal discretization scheme with the unit nodal director constraint enforced using the Lagrange multiplier method leads to a non-symmetric matrix 
while all the other three approaches lead to a symmetric matrix. 
We note that Table \ref{tab:matrix-forms} and this observation holds for Equations \eqref{eq:matrix_eq_iga}-\eqref{eq:matrix_eq_nodalR3S2-weak} which correspond to dynamics computations. 
In static cases, using the Lagrange multiplier method (\nodalSaddle), however, leads to a symmetric system matrix, 
while the system matrix resulting from the reduced equations using both Lagrange multiplier and nullspace methods (\nodalSaddleRed) remains non-symmetric.

In Figure \ref{fig:matrix_structure}, we illustrate the band structure 
of the system matrix on the left-hand side of the matrix Equations \eqref{eq:matrix_eq_iga}-\eqref{eq:matrix_eq_nodalR3S2-weak} for the five studied formulations. 
For this illustration, 
we consider an exemplary fixed-fixed cable commonly employed for airborne wind turbines, which is a cable DuPont's Kevlar 49 type 968. 
This type of cable has elastic constants of $E_{11} = 81.8$ GPa, 
a mass density of $\rho = 1.429 \cdot 10^3$ kg/m$^3$, 
an initial length of $L_0 = 300$ m, and 
a cross-sectional diameter of $0.007$ m.  
We discretize the cable with 20 elements in all cases and compute the matrices at the first time step and the second iteration to include the contribution of all terms. 
Since we want to focus on the symmetry and band structure of the system matrix, we employ cubic $C^1$ splines when using \iga \ to obtain the same number of dofs. 
For illustration purposes, we have removed the constrained degrees of freedom due to the fixed boundaries in all cases. 
Focusing on \iga \ (Figures \ref{fig:matrix_structure}a and b), 
we observe that 
the employed outlier removal approach \cite{hiemstra_outlier_2021} does not change the band structure or the symmetry of the system matrix, but only reduces its dimension since constraints are directly built into the spline space, as discussed in \cite{nguyen_rod_2024,hiemstra_outlier_2021}. 
Focusing on the nodal scheme with the unit nodal director constraint enforced using the Lagrange multiplier method (Figures \ref{fig:matrix_structure}c and d), 
we see the larger system matrix when using \nodalSaddle \ and a non-symmetric matrix when using \nodalSaddleRed. 
The matrix obtained with \nodalSaddle \ appears symmetric, however, is a non-symmetric matrix in dynamics computations since the top-right and bottom-left block matrices are evaluated at different time instances, as discussed in the previous section. 
Using the penalty method to enforce the unit nodal director constraint or neglecting it (Figures \ref{fig:matrix_structure}e and f, respectively) leads again to a symmetric matrix. 
Furthermore, we can see in Figure \ref{fig:matrix_structure} that using the isogeometric or nodal discretizations without the unit nodal director constraint enforced using the Lagrange multiplier method leads to a similar bandwidth of the system matrix. 
As discussed in \cite{nguyen_rod_2024}, this occurs regardless of the fact that spline basis functions employed for \iga \ have larger supports up to $p+1$ elements than that of nodal discretizations (see also \cite[p.~92-97]{Cottrell:09.1}). 
Notably, employing basis functions with larger support does not always lead to more evaluation per quadrature point. 
Using \iga, there are up to $p+1$ basis functions that have support in an element (see also Figure \ref{fig:spline_hermite}). 
For the studied rod formulation, using splines of $p<4$ leads to either less (if $p=2$) or the same (if $p=3$) number of active basis functions per element as using the nodal scheme.

In Table \ref{tab:matrix-forms}, we also highlight that using \iga\ with the strong outlier removal approach \cite{hiemstra_outlier_2021} or \nodalSaddleRed \ requires global matrix multiplication (see also Equations \eqref{eq:matrix_eq_iga_outlier_removal} and \eqref{eq:matrix_eq_nodalR3S2-strong-reduced}). 
This increases the computational cost, particularly for significantly large systems. 
On the one hand, using \iga \ with the strong outlier removal approach involves a constant multiplier, the extraction operator $\mat{C}$, which does not require any reassembly or update per iteration. 
Using \nodalSaddleRed, on the other hand, we need to reassemble the nullspace matrix $\mat{D}$ at each iteration since $\mat{D}$ depends on the nodal directors of the current configuration (see also Appendix \ref{sec:linearized_constraint_vector}). 
Moreover, we note that compared to \iga \ and \nodal, i.e. the nodal discretization scheme without considering unit nodal director constraint, the enforcement of this constraint requires the evaluation of additional terms, such as the constraint or the penalty terms, 
and hence requires the assembly of additional matrices and vectors on the left- and right-hand sides of the matrix equations.

\section{Numerical examples}\label{sec:results}
    
In this section, we investigate the accuracy and computational cost of the studied semi-discrete rod formulations discussed in the previous sections (see also Tables \ref{tab:semi-discrete-forms} and \ref{tab:matrix-forms}). 
To gain more insights, we first numerically study the condition number of the system matrix which depends on different parameters such as 
the penalty factor, the employed outlier removal approach, and 
the length of the nodal directors. 
We then numerically illustrate via an example of a planar roll-up that preserving nodal directors in the unit sphere leads to better accuracy than nodal directors in $\mathbb{R}^3$. 
Via a static and dynamic analysis of exemplary cables, 
we show that cubic $C^1$ isogeometric discretization leads to the same responses in the static case, however, slightly larger responses in parts of the dynamic computation. 
We also illustrate that 
\iga \ with or without outlier removal averagely requires a smaller computational time per iteration than using any formulations based on the nodal scheme. 
For fine meshes, using \nodalSaddleRed \ requires the most time per iteration compared to other approaches. 
We also numerically demonstrate via these examples 
that using the nodal scheme with unit nodal director constraint enforced using the Lagrange multiplier method leads to zero nodal axial stress resultants, as discussed in the previous sections.

\subsection{Numerical study of the condition number of the system matrix}\label{sec:cond_number}

The condition number of the system matrix in the matrix equations plays an essential role in ensuring the convergence of the Newton-Raphson iterative procedure, which is associated with the robustness of the corresponding formulation. 
Hence, 
to gain better insights into the influence of different parameters such as the penalty factor, the outlier removal approach, and the length of the nodal directors on this number, as well as insights into the robustness of each formulation for the same number of elements, 
we numerically study and compare this obtained with five semi-discrete formulations considered in this work (see also Tables \ref{tab:semi-discrete-forms} and \ref{tab:matrix-forms}). 
We consider an exemplary cable commonly employed for airborne wind turbines, which is a cable DuPont's Kevlar 49 type 968. 
This type of cable has elastic constants of $E_{11} = 81.8$ GPa,  
a mass density of $\rho = 1.429 \cdot 10^3$ kg/m$^3$, 
an initial length of $L_0 = 300$ m, and 
a cross-sectional diameter of $0.007$ m. 
The cable is fixed on both ends and subjected to its self-weight. 
We discretize the cable with 40 elements using five formulations discussed in the previous sections. 
Since we focus on the condition number, we choose cubic $C^1$ splines for the isogeometric discretization to obtain the same number of dofs as the nodal scheme. 
We compute the condition number of the system matrix at the first load step and the first iteration for all formulations.

\begin{figure}[ht]
	\centering
    \def\svgwidth{1\textwidth}
    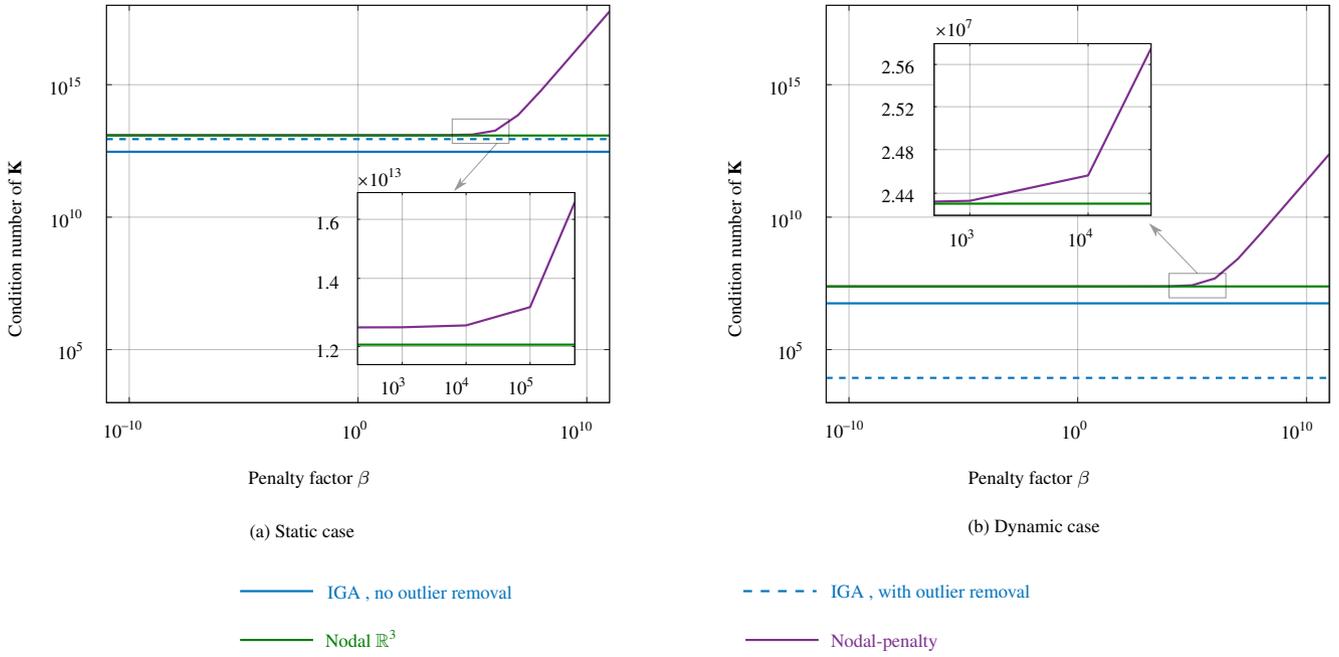
 \caption{Condition number of the system matrix as a function of the penalty factor $\beta$.}\label{fig:condition_number_penalty}
\end{figure}

\begin{figure}[ht]
	\centering
    \def\svgwidth{1\textwidth}
    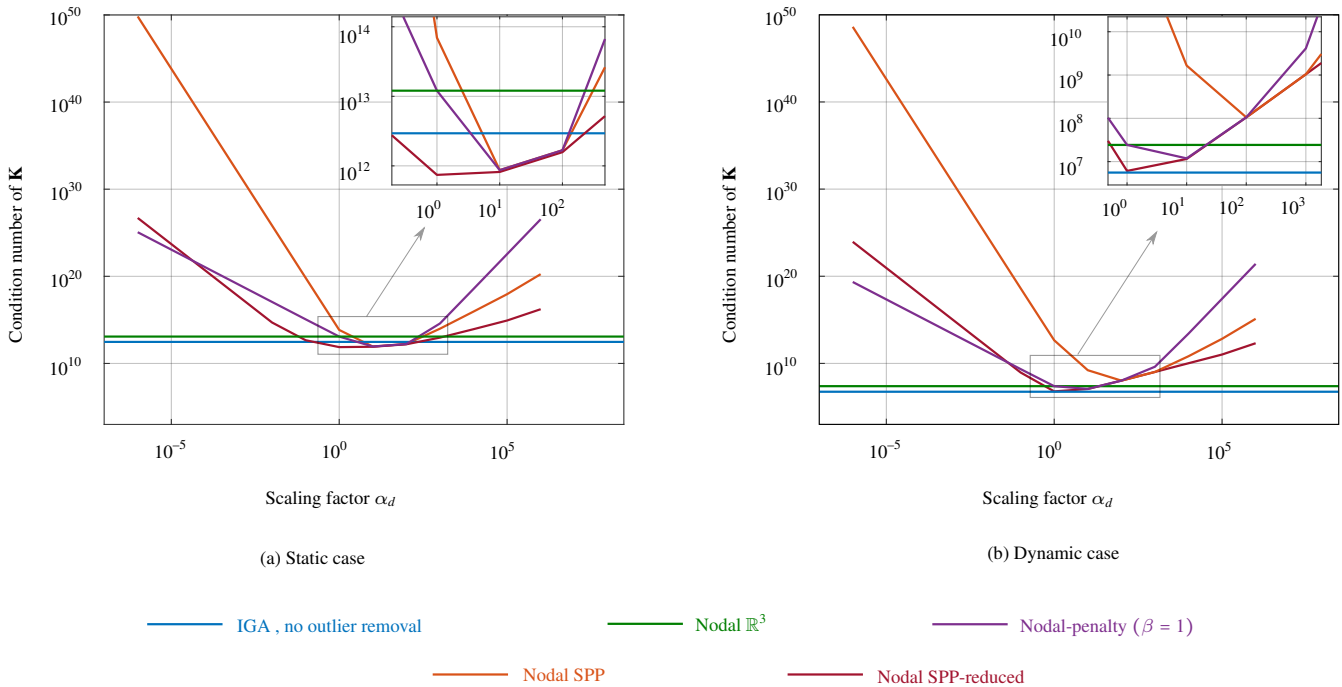
 \caption{Condition number of the system matrix as a function of the enforced length for the nodal directors.}\label{fig:condition_number_scalingD}
\end{figure}

We first investigate the effect of the penalty factor on the condition number of the system matrix when using \nodalPenalty. 
In Figure \ref{fig:condition_number_penalty}, we plot this number as a function of the penalty factor $\beta$ in a logarithmic scale. 
We compare the condition number obtained with \nodalPenalty \ (purple curve), \nodal \ (green line), and \iga \ with (blue dashed line) and without (blue continuous line) outlier removal. 
We consider the system matrix employed for both static (Figure \ref{fig:condition_number_penalty}a) and dynamic cases (Figure \ref{fig:condition_number_penalty}b). 
Focusing on the purple curve, we observe that increasing the penalty factor $\beta$ increases the condition number in both cases. 
This necessarily means that a significantly large penalty factor can lead to ill-conditioning of the system matrix, reducing the robustness of the corresponding formulation \nodalPenalty. 
For the studied cable, the penalty terms start having an effect on the condition number of the system matrix when the penalty factor is larger than $10^6$ (see also inset figures in Figure \ref{fig:condition_number_penalty}). 
Using smaller penalty factors leads to approximately the same condition number as using \nodal \ , as expected. 
Focusing on the green and blue lines, 
we see that using \iga \ generally leads to a smaller condition number than \nodal. 
In the static case, employing the strong approach of outlier removal \cite{hiemstra_outlier_2021} slightly increases this to approximately the same condition number as using \nodal. 
In the dynamic case, however, the outlier removal approach reduces the condition number by several orders of magnitude. 
This is consistent with the observations and discussions in \cite{nguyen_rod_2024} that using the outlier removal approach \cite{hiemstra_outlier_2021} improves the robustness of isogeometric discretizations for dynamics computations. 
Based on this observation, in this work, we generally employ the outlier removal approach only for dynamic computations, unless it is stated otherwise. 
Furthermore, comparing the condition number obtained in the static and dynamic cases, we can see that 
the obtained condition number in the dynamic case is several orders of magnitude smaller than that in the static case.  
This necessarily means that the terms regarding the mass matrix improve the conditioning of the system matrix, improving the robustness of the studied formulations. 
We conclude that the penalty factor should be chosen as large as required to enforce the unit nodal director constraint such that it does not negatively affect the condition number of the resulting system matrix. 
For dynamics computations, the outlier removal approach \cite{hiemstra_outlier_2021} improves the conditioning of the system matrix, leading to better robustness when using \iga, as discussed in \cite{nguyen_rod_2024}.

A well-established approach 
to improve the conditioning of the system matrix for thin-walled structures is the scaled director conditioning \cite{Kloppel_scaled_director2011,Wall_scaled_director2000,Gee_scaled_director2005}. 
The idea is to scale the length of the director variable field to reduce the condition number of the resulting system matrix. 
Hence, we also investigate whether this approach improves the conditioning for the studied semi-discrete formulations when using the nodal discretization scheme. 
We employ this approach by replacing the constrained nodal directors $\vect{d}_h$ in the discretization \eqref{eq-discretize_nodal} by $\frac{1}{\alpha_d}\hat{\vect{d}}_h$ and enforcing a length of $\alpha_d$ for $\hat{\vect{d}}_h$. 
We note that the resulting nodal directors are then still constrained to have a unit length. 
In Figure \ref{fig:condition_number_scalingD}, 
we plot the condition number of the system matrix as a function of the scaling factor $\alpha_d$ when using \nodalPenalty \ (purple curve), \nodalSaddle \ (orange curve), and \nodalSaddleRed \ (dark red curve). 
When using \nodalPenalty, to avoid the effect of the penalty factor on the condition number, we choose a penalty factor $\beta$ such that $\beta \,(2EI/L_0) = 1.0$. 
For comparison purposes, we also include the condition number obtained when using \nodal \ (green line) and \iga \ without outlier removal (blue line). 
The inset figures focus on the region close to the minimum of the condition number with respect to $\alpha_d$.

Focusing on the purple curve obtained with \nodalPenalty, we observe that we obtain a minimal value with $\alpha_d = 10$. 
Increasing $\alpha_d$ increases the condition number 
that is then several orders of magnitude larger than that obtained with \nodalSaddle \ or \nodalSaddleRed \ using the same value of $\alpha_d > 10$. 
Focusing on the dark red curve obtained with \nodalSaddleRed, we see that we obtain a minimal value without scaling the enforced length, i.e. $\alpha_d = 1.0$. 
Compared to the minimum obtained with \nodalPenalty \ or \nodalSaddle, the obtained value is the smallest minimum that is slightly smaller than the condition number obtained with \iga \ in the static case and is the same as that in the dynamic case. 
Comparing the dark red and orange curves, we also see that reducing the system of equations using the nullspace method significantly improves the conditioning of the system matrix. 
Focusing on the orange curve obtained with \nodalSaddle, 
we observe that we obtain a minimal condition number with $\alpha_d = 10$ and $\alpha_d = 100$ in the static and dynamic cases, respectively. 
The obtained minimal value is the same as the minimum when using \nodalPenalty \ in the static case and is one order of magnitude larger in the dynamic case. 
We conclude that when enforcing the unit nodal director constraint using the Lagrange multiplier method, reducing the system of equations using nullspace method (\nodalSaddleRed) also improves the conditioning of the system matrix and hence the robustness of the corresponding formulation. 
One does not need to employ 
the scaled director conditioning approach to further reduce the condition number. 
Using \nodalSaddleRed \ can lead to a similar condition number as using \iga. 
When enforcing the unit nodal director constraint using \nodalSaddle \ or \nodalPenalty, the scaled director conditioning approach can improve the conditioning, however, requires a parameter study in advance to estimate optimal scaling factors.

\subsection{Convergence study of a planar roll-up}\label{sec:convergence_study_rollup}

\begin{figure}[ht]
	\centering
    \def\svgwidth{0.5\textwidth}
    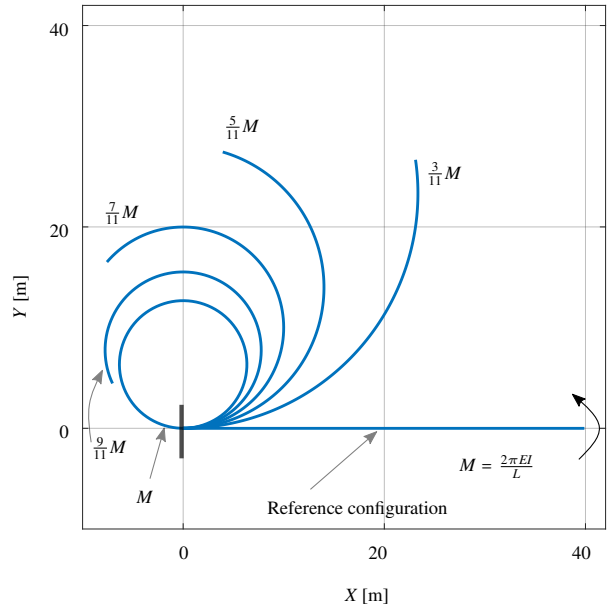

 \caption{Deformed configurations of the clamped rod bent to a circle at different load steps, computed with \textbf{cubic} $C^1$ isogeometric discretization without outlier removal approach and a mesh of 8 elements.} \label{fig:snapshots_rollup}
\end{figure}

\begin{figure}[ht]
	\centering
    \def\svgwidth{0.9\textwidth}
    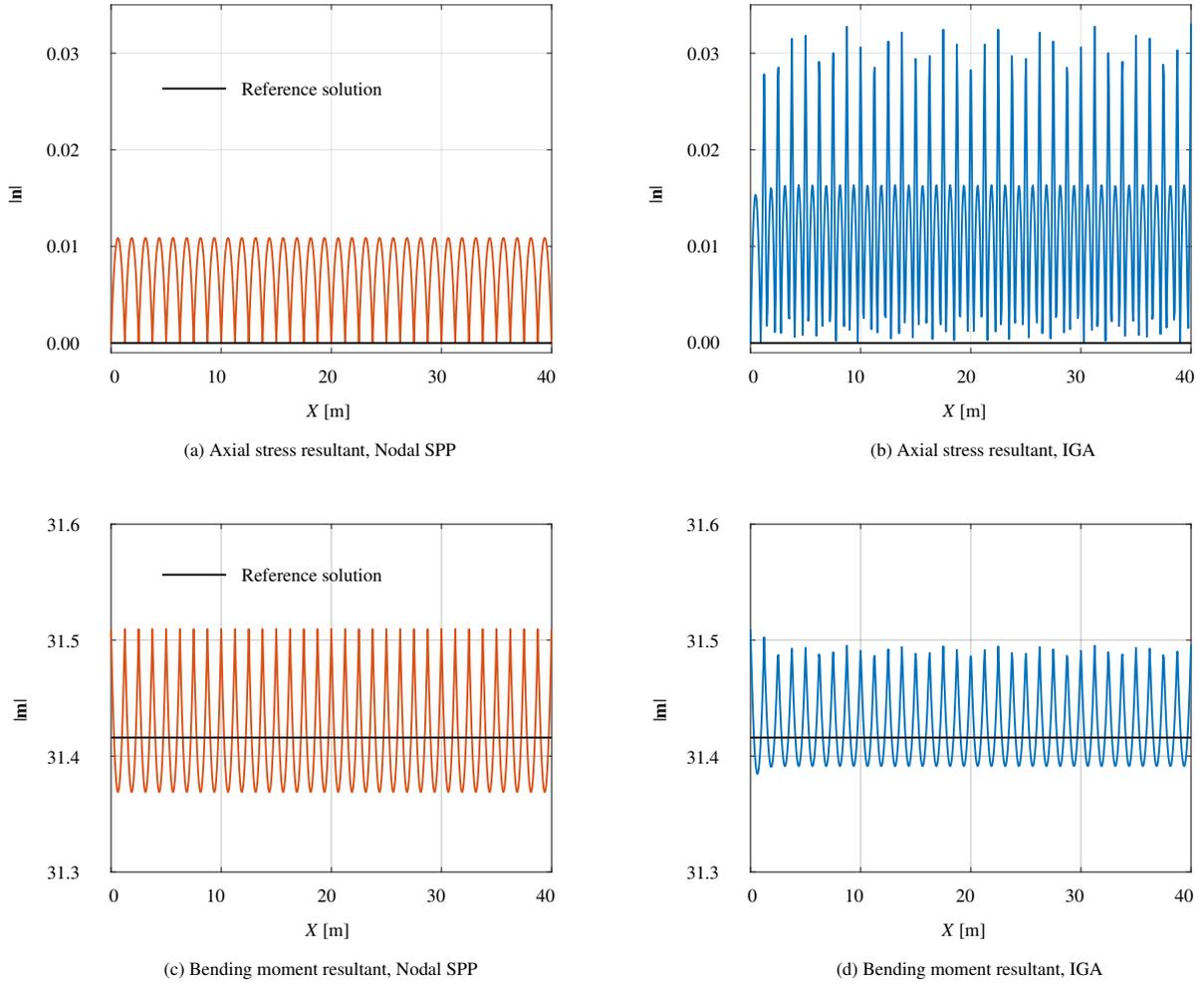

 \caption{Axial stress and bending moment resultants in the final configuration of a planar roll-up, computed with \textbf{cubic} isogeometric and nodal discretizations.} \label{fig:stress_compare_rollup}
\end{figure}

\begin{figure}[ht]
	\centering
    \def\svgwidth{0.9\textwidth}
    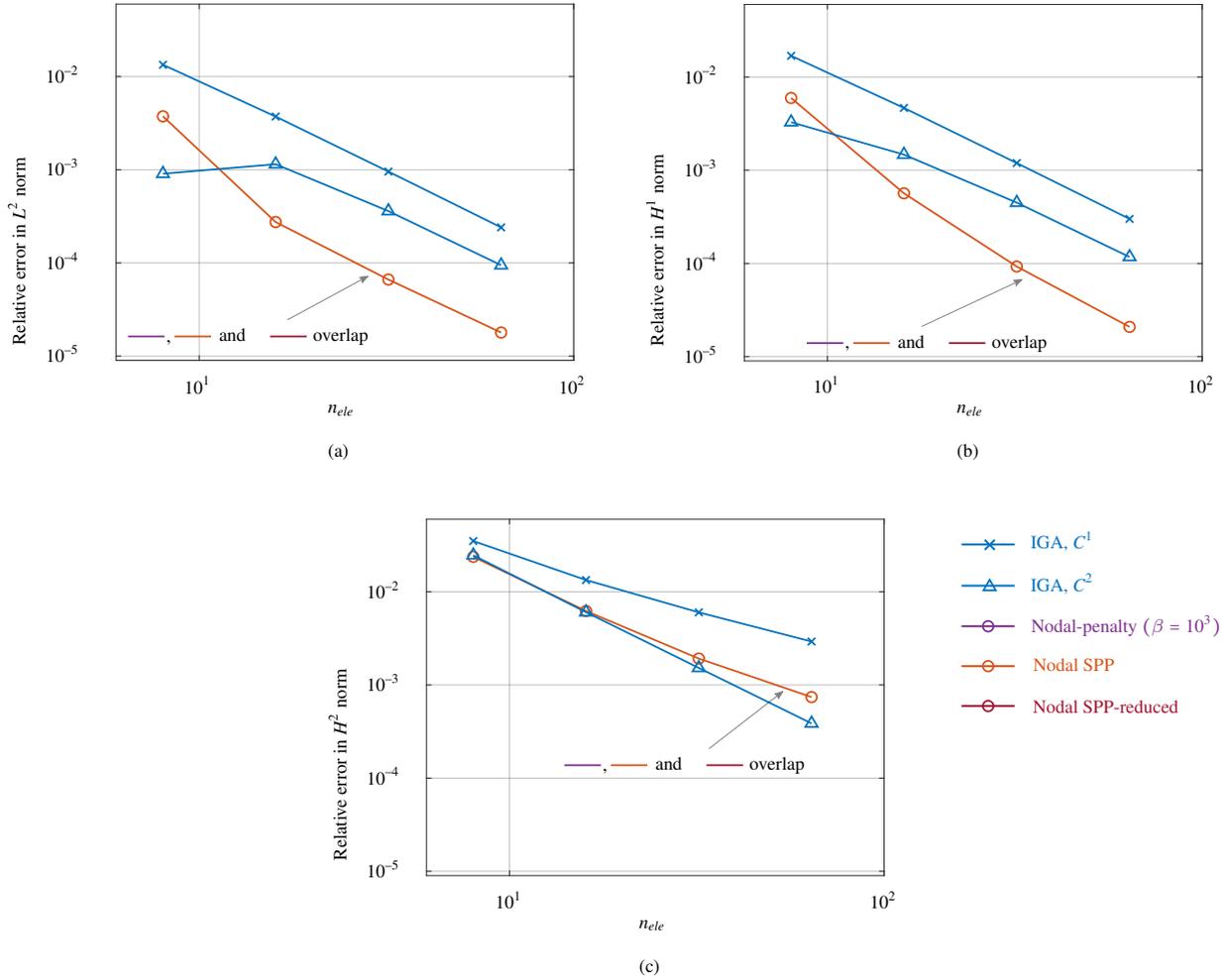

 \caption{Convergence of relative errors of the clamped rod bent to a circle computed with \textbf{cubic} isogeometric and nodal discretizations.} \label{fig:convergence_compare_rollup}
\end{figure}

We now investigate the accuracy and convergence achieved with the studied semi-discrete formulations. 
To this end, we consider 
a planar roll-up example of pure bending with a well-known final configuration and stress resultants. 
Since we focus on the accuracy achieved with different spatial discretization schemes, we employ cubic $C^1$ splines when using \iga \ throughout this subsection, which belong to the same space as the Hermit splines employed for nodal discretization schemes. 
We consider a rod with 
an initital length of $L=40$ m, 
axial stiffness of $EA=100$ N, and 
bending stiffness of $EI=200$ Nm$^2$, which we discretize with 40 elements. 
The rod is clamped at the left-end and subjected to a bending moment $M = 2 \pi EI / L$ that is to roll the rod to one full circle and is applied in 55 load steps, illustrated in Figure \ref{fig:snapshots_rollup}. 
We strongly enforce the Dirichlet boundary conditions in the standard way and 
employ a tolerance of $10^{-10}$ for the convergence of the Newton-Raphson method in all cases. 
Figure \ref{fig:snapshots_rollup} also shows five snapshots of the deformed configurations, computed with the cubic $C^1$ isogeometric discretization. 
We obtained the same configurations when using \nodalPenalty, \nodalSaddle, or \nodalSaddleRed, 
and hence discarded these in Figure \ref{fig:snapshots_rollup} purely for the sake of illustration clarity. 
For this example, using \nodal \ leads to an ill-conditioned system matrix and thus the Newton-Raphson scheme does not converge. 
When using \nodalPenalty, based on empirical results, we observed that 
increasing the penalty factor $\beta < 10^9$ improves the conditioning but  
using small values of $\beta \leq 10^1$ does not sufficiently improve the conditioning to achieve convergence of the Newton-Raphson scheme. 
We also saw that using $\beta \leq 10^2$ or $\beta \geq 10^7$ requires more than 10 iterations per load step. 
Hence, we choose a penalty factor of $\beta=10^5$ for our computations, which requires a maximum of 8 iterations per load step. 
Using \iga, \nodalSaddle, or \nodalSaddleRed \ requires a maximum of 6 iterations per load step.

Figure 
\ref{fig:stress_compare_rollup} illustrates the axial stress resultants (Figures 
\ref{fig:stress_compare_rollup}a and b) and the bending moment resultants (Figures 
\ref{fig:stress_compare_rollup}c and d) 
in the circular configuration at the last load step, 
obtained with \iga \ and \nodalSaddle. 
Using \nodalSaddleRed \ or \nodalPenalty \ with $\beta=10^5$ leads to the same stress resultants as \nodalSaddle. 
Hence, we compare the results obtained with \iga \ and \nodalSaddle \ in Figure 
\ref{fig:stress_compare_rollup}. 
We also include the reference value (black curves) of the stress resultants for this pure bending example. 
Focusing on the bending moment resultants in Figures 
\ref{fig:stress_compare_rollup}c and d, 
we observe that both discretization schemes lead to stress resultants consisting of oscillations around the reference value along the rod. 
Focusing on the axial stress resultants in Figures \ref{fig:stress_compare_rollup}a and b, 
we see that both schemes lead to non-zero and oscillating stress resultants in this example. 
When using \nodalSaddle, due to the enforced unit nodal director constraint, 
we do expect oscillations in non-zero axial stress resultants since the nodal stress values are constrained to zero, as discussed in Section \ref{sec:discrete_rod_space}. 
When using \iga, we observe the same results despite the absence of the unit nodal director constraint (see Figure \ref{fig:stress_compare_rollup}b). 
We note that based on empirical results, increasing the polynomial degree or continuity order when using \iga, or refining the mesh when using either one of the four aforementioned formulations reduces the magnitude of the oscillations in the stress resultants. 
These results imply the typical effect of membrane locking on stress resultants \cite{Stolarski1982}. 
We note that for 
geometrically nonlinear problems, another locking phenomenon has been reported for 
shell structures in \cite{Willmann_nonlinear_locking2023}. 
This requires further study in future work to identify whether this occurs for the studied rod formulation.

Figure \ref{fig:convergence_compare_rollup} illustrates the convergence of the relative errors in $L^2$-, $H^1$-, and $H^2$-norms of the studied planar roll-up, obtained with isogeometric and nodal discretizations via uniform $h$-refinement. 
We again obtain the same accuracy and convergence behavior in all three error norms when using \nodalPenalty, \nodalSaddle, and \nodalSaddleRed. 
Hence, we plot the results obtained with \nodalSaddle \ (orange curves with circle markers) and compare with those obtained with \iga \ based on cubic $C^1$ (blue curves with cross markers) and $C^2$ (blue curves with triangle markers) splines. 
Focusing on the blue curves, we observe that increasing the continuity order of cubic spline functions increases the accuracy in all cases, as discussed also in \cite{nguyen_rod_2024}. 
Focusing on the convergence in the $L^2$- and $H^1$-norms (Figures \ref{fig:convergence_compare_rollup}a and b), we see that using \iga \ 
with cubic $C^1$ splines, on the one hand, yields a convergence rate similar to that of \nodalSaddle, on the other hand, leads to a higher error level. 
In this example, one possible reason for the offset between these two convergence curves might be the enforced unit nodal director constraint. 
Since the reference solutions include zero axial stress resultants, enforcing the unit nodal director constraint when using, e.g., \nodalSaddle\ ensures the solutions at the nodes, which is not the case when using \iga, and thus leads to smaller errors. 
Increasing the continuity of the spline basis functions reduces the error level, approximately to the same level as \nodalSaddle. 
For the errors in the $H^2$-norm (Figure \ref{fig:convergence_compare_rollup}c), we observe that using \iga\ with cubic $C^2$ splines leads to a lower error level of almost one order of magnitude and a slightly higher convergence rate. 
We conclude that using the nodal discretization scheme that preserves the nodal director field in the unit sphere  
leads to better accuracy in the deformations in different error norms than using cubic $C^1$ \iga. Increasing the continuity of the employed cubic spline basis functions further increases both accuracy and convergence rate. 
The nodal scheme without this property, i.e. \nodal \ or \nodalPenalty \ with small penalty factors that do not sufficiently improve the conditioning, may lead to an ill-conditioned system matrix and hence is generally less robust than \iga, \nodalSaddle, and \nodalSaddleRed.

\begin{remark}
    In our previous work \cite{nguyen_rod_2024}, we have studied the convergence and accuracy of isogeometric discretizations with increasing polynomial degree and continuity order. For additional results and further discussions, we refer to \cite{nguyen_rod_2024}.
\end{remark}

Furthermore, we also study the convergence of the considered planar roll-up with different slenderness ratios in Appendix \ref{sec:convergence_rollup_radi}. 
We observe the typical pre-asymptotic plateau in the convergence curves for large slenderness ratios, which becomes more severe with increasing slenderness ratios in all cases. 
This implies the effect of membrane locking on the accuracy of the spatial discretization schemes, which is consistent with the observations in the stress resultants discussed above. 
In this work, we focus on the comparison between different semi-discrete rod formulations and hence plan to investigate this locking effect and locking-preventing techniques for the studied rod formulation in future work.

\subsection{Static analysis of a catenary}

\begin{figure}[ht]
	\centering
    \def\svgwidth{0.55\textwidth}
    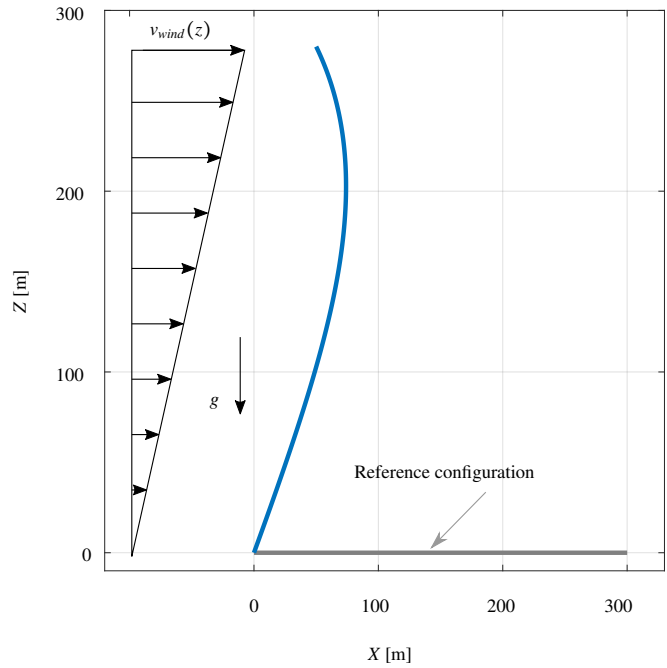

 \caption{Final deformed configuration of an exemplary cable subjected to its weights and a linear wind profile with a prescribed fairlead position, computed with cubic $C^1$ isogeometric discretization without outlier removal.} \label{fig:final_config_compare_catenary}
\end{figure}

\begin{figure}[ht]
	\centering
    \def\svgwidth{0.9\textwidth}
    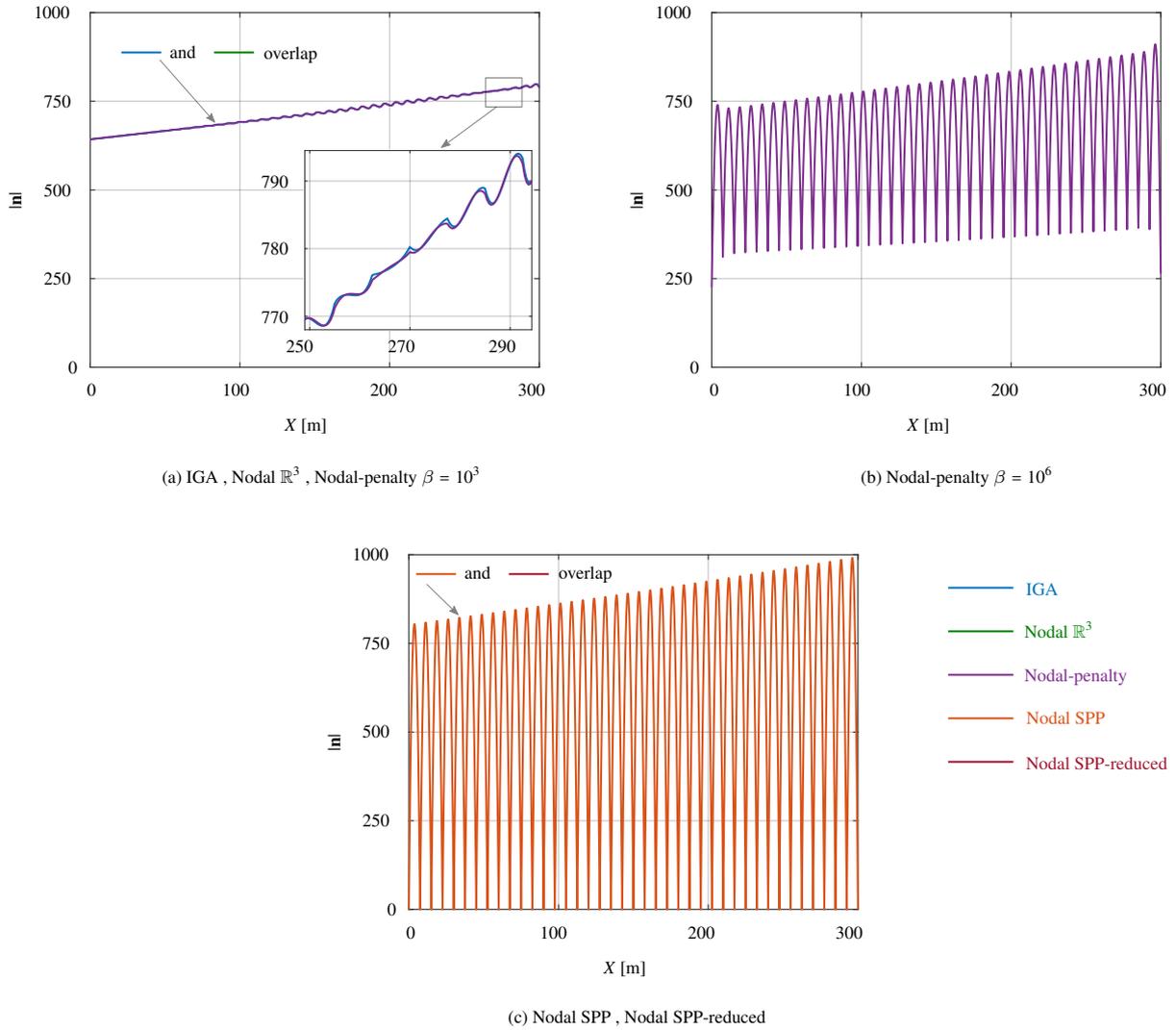
	
	\caption{Axial stress resultants in the final configuration of the cable in Figure \ref{fig:final_config_compare_catenary}, computed with isogeometric and nodal discretization schemes.} \label{fig:membrane_compare_catenary}
\end{figure}

\begin{figure}[ht]
	\centering
    \def\svgwidth{0.6\textwidth}
    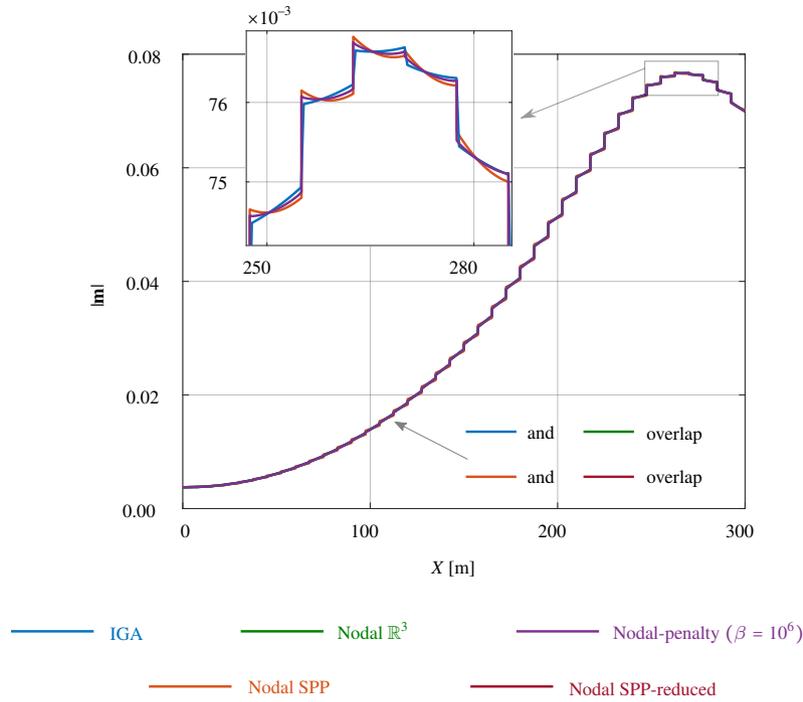    
	
	\caption{Bending moment resultants in the final configuration of the cable in Figure \ref{fig:final_config_compare_catenary}, computed with isogeometric and nodal discretization schemes.} \label{fig:bending_compare_catenary}
\end{figure}

\begin{table}[ht]
    \centering
    \begin{tabularx}{1\linewidth}{|>{\centering\arraybackslash}X | >{\centering\arraybackslash}X | >{\centering\arraybackslash}X | >{\centering\arraybackslash}X | >{\centering\arraybackslash}X | >{\centering\arraybackslash}X |}
        \toprule
        Number of elements & \iga$^1$ & \nodal & \nodalSaddle & \nodalSaddleRed & \nodalPenalty$^2$ \\
        \hline
        $2^3$ & 21 & 21  & 12 & 24 & 20 \\
        $2^4$ & 13 & 13  & 13 & 13 & 13 \\
        $2^5$ & 21 & 21  & 17 & 15 & 13 \\
        $2^6$ & 18 & $-$$^3$ & 18 & 18 & 20 \\
        $2^7$ & 17 & 17  & 26 & 22 & 31 \\
        $2^8$ & 23 & 23  & 39 & 25 & 21 \\
        \hline
        \multicolumn{6}{|l|}{\textit{$^1$ The maximum number of iterations is the same for both cases with and without outlier removal.}} \\
        \multicolumn{6}{|l|}{\textit{$^2$ $\beta=10^8$.}} \\
        \multicolumn{6}{|l|}{\textit{$^3$ The system matrix is ill-conditioned. The Newton-Raphson scheme did not converge.}} \\
        \bottomrule
    \end{tabularx}
     \caption{Maximum number of iterations for the Newton-Raphson scheme when using different discretization schemes. This occurred once at the 52$^{nd}$ load step, when we start moving the right-end of the cable in Figure \ref{fig:final_config_compare_catenary} to the final fairlead position.}\label{tab:static_max_iter}
\end{table}

\begin{figure}[ht]
	\centering
    \def\svgwidth{0.85\textwidth}
    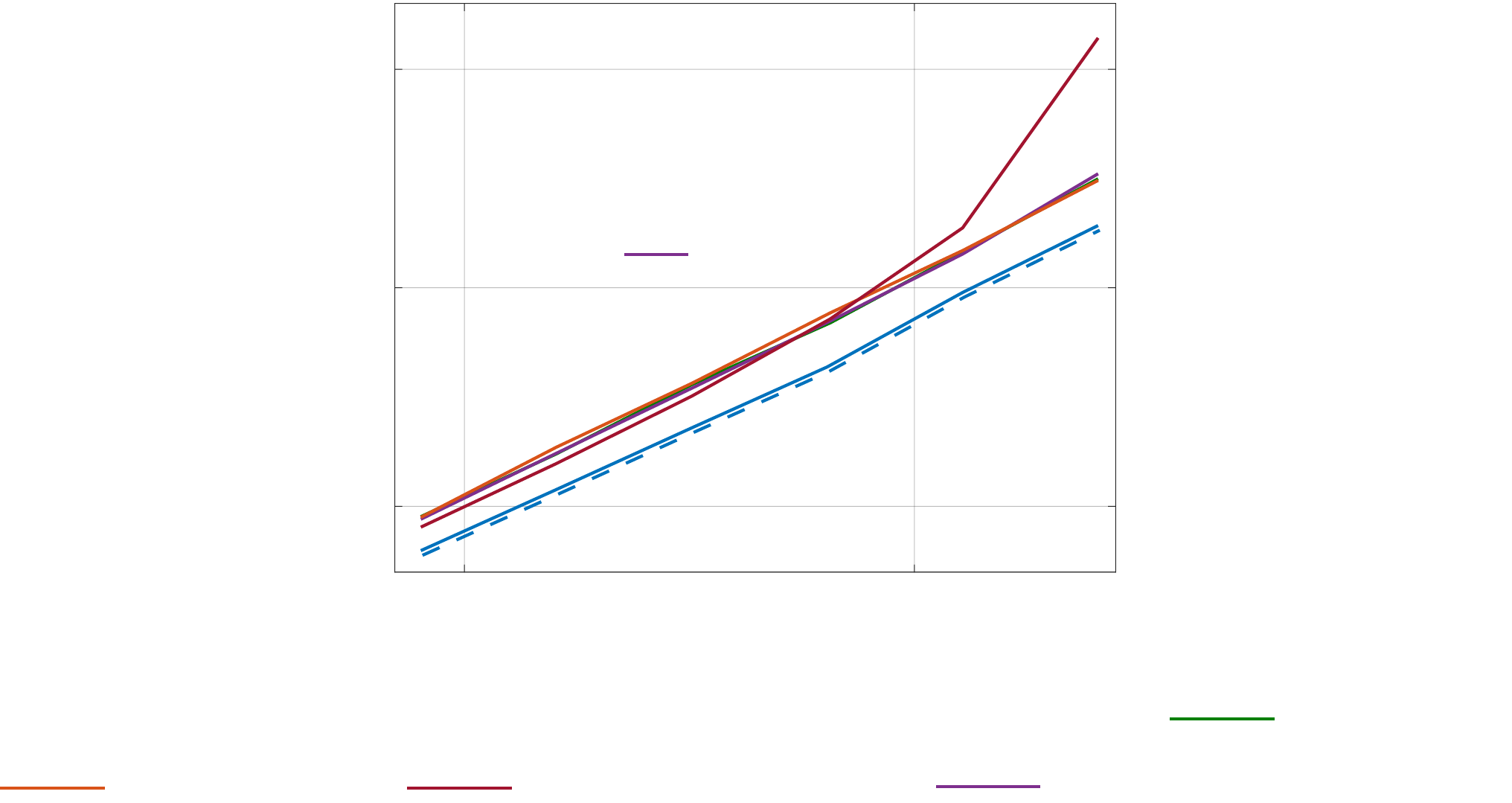    
	
	\caption{Computation time per iteration, averaged over all iterations and load steps, required when using isogeometric and nodal discretization schemes for the cable in Figure \ref{fig:final_config_compare_catenary}.}\label{fig:static_averaged_computing_time}
\end{figure}

We now investigate the performance of the considered semi-discrete rod formulations, in terms of the obtained responses and computational time, for examples that are relevant and common for studying mooring lines in offshore wind engineering. 
The first example is a static analysis of an exemplary cable of the DuPont's Kevlar 49 type 968, which is also considered in Section \ref{sec:cond_number}. 
The cable lies initially straight along the $X$-axis and is fixed at its left-end while we move the right-end 
from the initial position of $(L_0,0,0)$ to 
a prescribed fairlead position of $(50,0,280)$. 
We consider the self-weight of the cable and the following linear wind profile along the $Z$-axis (see also Figure \ref{fig:final_config_compare_catenary}):
\begin{align}
    v_{\text{wind}} (z) = 15 \, \frac{z}{100} \; [\text{m/s}] \,.
\end{align} 
To activate the stiffness of the cable, 
we first prestress the cable by moving the right-end to $(L_0+0.01,0,0)$ in one load step. 
We then apply the self-weight in 50 load steps while keeping the right-end at the prestressed position. 
After the self-weight is fully applied, we move the right-end to the prescribed fairlead position in 400 load steps. 
We compute the response of this cable using the five semi-discrete formulations: \iga, \nodal, \nodalPenalty, \nodalSaddle, and \nodalSaddleRed. 
We discretize the cable in 40 elements and employ 
a tolerance of $10^{-10}$ for the convergence of the Newton-Raphson method in all cases. 
We again focus on the performance of different discretization schemes and hence employ basis functions in the same function space, i.e. cubic $C^1$ and Hermite splines. 
When using nodal discretizations, we enforce the homogeneous boundary condition at the left-end in the standard way but using the extraction operator to enforce this when using \iga \ for the sake of implementation when employing the outlier removal approach \cite{hiemstra_outlier_2021}. 
Nevertheless, the prescribed right-end is enforced at each step in the standard way for all formulations. 
The maximum number of iterations for all formulations is 18 when we start moving the right-end to the fairlead position, i.e. at the 52$^{nd}$ load step. 
For all remaining load steps, all formulations require a maximum of 6 iterations.

Figure \ref{fig:final_config_compare_catenary} illustrates the final configuration of the studied cable obtained with cubic $C^1$ isogeometric discretization without outlier removal. 
Using either \nodal, \nodalPenalty, \nodalSaddle, \nodalSaddleRed, or \iga \ with outlier removal leads to the same configuration. 
Hence, for illustration clarity, we plot only one result in Figure \ref{fig:final_config_compare_catenary}. 
In Figure \ref{fig:membrane_compare_catenary}, we plot the axial stress resultants obtained with each formulation. 
When using \iga, the outlier removal approach does not affect the accuracy (see also \cite{nguyen_rod_2024,hiemstra_outlier_2021}), we refer to \iga \ for both cases with and without outlier removal in the following. 
We also include an overkill solution (gray curve) computed with \iga\ using cubic $C^1$ splines and 2048 elements as a reference solution. 
Focusing on the blue and green curves obtained with \iga \ and \nodal, respectively, in Figure \ref{fig:membrane_compare_catenary}a, 
we observe that these two schemes lead to the same result. 
This is expected since we employ basis functions in the same function space for these schemes. 
We note that the axial stress resultants obtained with these two formulations consist of slight oscillations 
which can be eliminated by refining the mesh (see the gray curve in Figure \ref{fig:membrane_compare_catenary}a). 
The cause of these oscillations might be membrane locking, as discussed in the previous example of a planar roll-up. 
We also see that using \nodalPenalty \ with a penalty factor of $\beta=10^5$ (purple curve in Figure \ref{fig:membrane_compare_catenary}a), the unit nodal director constraint is not enforced sufficiently, leading approximately to the same axial stress resultants as using \iga \ or \nodal. 
Increasing $\beta$ enforces this constraint more sufficiently, however, leads to constrained nodal axial stress values and thus to oscillations in the axial stress resultants, as illustrated in Figure \ref{fig:membrane_compare_catenary}b for a larger value of $\beta=10^8$. 
We note that using 
a penalty factor of $\beta=10^8$, nevertheless, does not sufficiently enforce the unit nodal director constraint and hence the nodal axial stresses are not zero. 
Enforcing this constraint using either \nodalSaddle \ or \nodalSaddleRed \ leads to axial stress resultants with zero nodal values, i.e. oscillating resultants, as discussed in Section \ref{sec:discrete_rod_space}. 
We also see that reducing the system of equations of \nodalSaddle \ using the nullspace method does not affect the accuracy and thus leads to the same results. 
Figure \ref{fig:bending_compare_catenary} illustrates the bending moment resultants obtained with the five aforementioned formulations. 
We observe that all formulations lead to approximately the same result. 
We see a slight difference in nodal values due to the enforced unit nodal director constraint when using either the Lagrange multiplier or penalty method. 
We also observe that refinining the mesh eliminates oscillations in the bending moment resultants (see the gray curve in Figure \ref{fig:bending_compare_catenary}), which might result from the effect of membrane locking.

In terms of computational cost, we investigate the maximal number of iterations and computing time per iteration when using the studied formulations under mesh refinement. 
In Table \ref{tab:static_max_iter}, we list the maximum number of iterations required for the 52$^{nd}$ load step. 
Focusing on the first two columns, we observe that using \iga \ and \nodal \ requires the same number of iterations. 
Using \nodal, however, can lead to an ill-conditioned system matrix, e.g. for the mesh of $2^6$ elements for the studied cable. 
Compared with these two formulations, 
using \nodalSaddle \ requires fewer iterations for coarser meshes and, nonetheless, more iterations for finer meshes. 
We also see that using \nodalSaddleRed \ and \nodalPenalty \ requires approximately the same number of iterations as \iga \ or \nodal. 
In Figure \ref{fig:static_averaged_computing_time}, we plot the averaged computing time per iteration and load step when using the aforementioned formulations. 
We observe that using \iga \ (blue curves) requires slightly less time for computing per iteration than using any of the nodal schemes. 
We see that using outlier removal (blue dashed curve) requires approximately the same time as without outlier removal (blue solid curve). 
This is expected since we employ the extraction operator to enforce the homogeneous boundary condition, i.e. we perform a global matrix multiplication in each iteration in both cases. 
We note that although we enforce the boundary condition in a more expensive way when using \iga \ than the standard way when using the nodal scheme, \iga \ requires less time per iteration than the nodal scheme. 
This necessarily means that for this example, the global multiplication is not decisive for the computational cost per iteration. 
Focusing on the green, orange, and purple curves, 
we observe that using \nodal, \nodalSaddle, and \nodalPenalty \ requires the same computing time per iteration. 
Using \nodalSaddleRed \ requires approximately the same time on coarser meshes, however, significantly more computational effort on finer meshes. 
This is due to the reassembly of the nullspace matrix $\mat{D}$ at each iteration, as discussed in Section \ref{sec:computational_cost}.

We conclude that for the static analysis of the studied cable, the cubic isogeometric and nodal discretizations show approximately the same accuracy in the deformed configuration and bending moment resultants. 
A nodal scheme with unit nodal director constraint enforced using the Lagrange multiplier method leads to zero nodal axial stresses and hence oscillating axial stress resultants. 
A nodal scheme without this constraint leads to the same axial stress as the isogeometric scheme. 
In general, all studied semi-discrete formulations require a similar number of iterations for this example. 
However, \nodal \ seems to be the least robust formulation compared to others since it can lead to an ill-conditioned system matrix. 
Per iteration, \iga \ requires the smallest computing time despite the global matrix multiplication with the constant extraction operator. 
\nodal, \nodalPenalty, and \nodalSaddle \ require approximately the same computing time despite the evaluation of additional penalty terms and a larger system of equations of the saddle-point problem. 
\nodalSaddleRed \ requires approximately the same computing time as these three formulations on coarser meshes, however, significantly more computational effort on finer meshes.

\begin{figure}[ht]
	\centering
    \def\svgwidth{0.7\textwidth}
    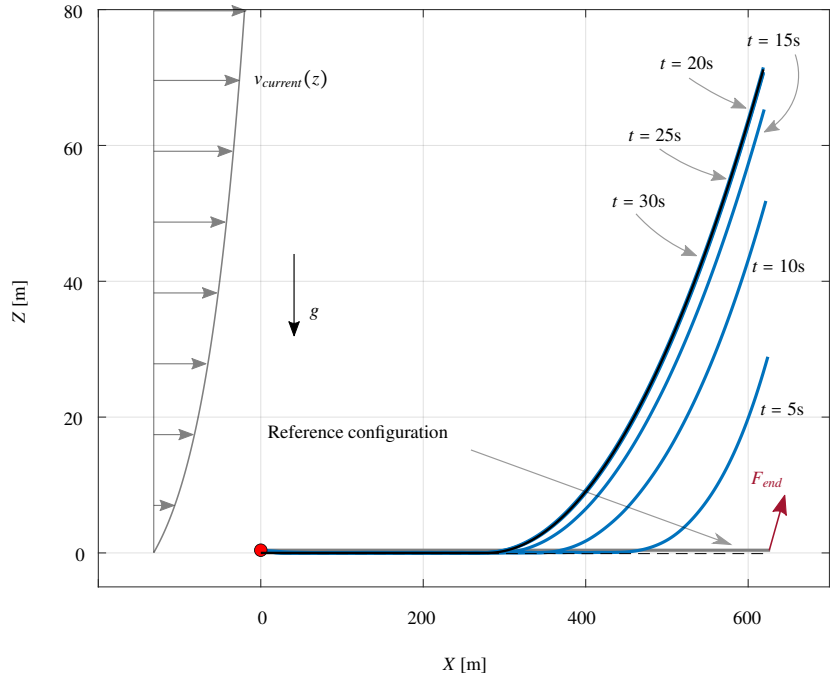

 \caption{Snapshots of a mooring line subjected to a logarithmic current speed profile in open water, computed with cubic $C^1$ isogeometric discretization and the outlier removal approach \cite{hiemstra_outlier_2021}. The snapshots in the last 10 seconds and the final configuration from \cite[Ch.~7.6.4, p.~257]{Gunnar2024} (black solid curve) are not distinguishable in this axis scale. The horizontal black dashed line represents the numerical barrier.} \label{fig:final_config_compare_mooringline}
\end{figure}

\begin{figure}[ht]
	\centering
    \def\svgwidth{0.87\textwidth}
    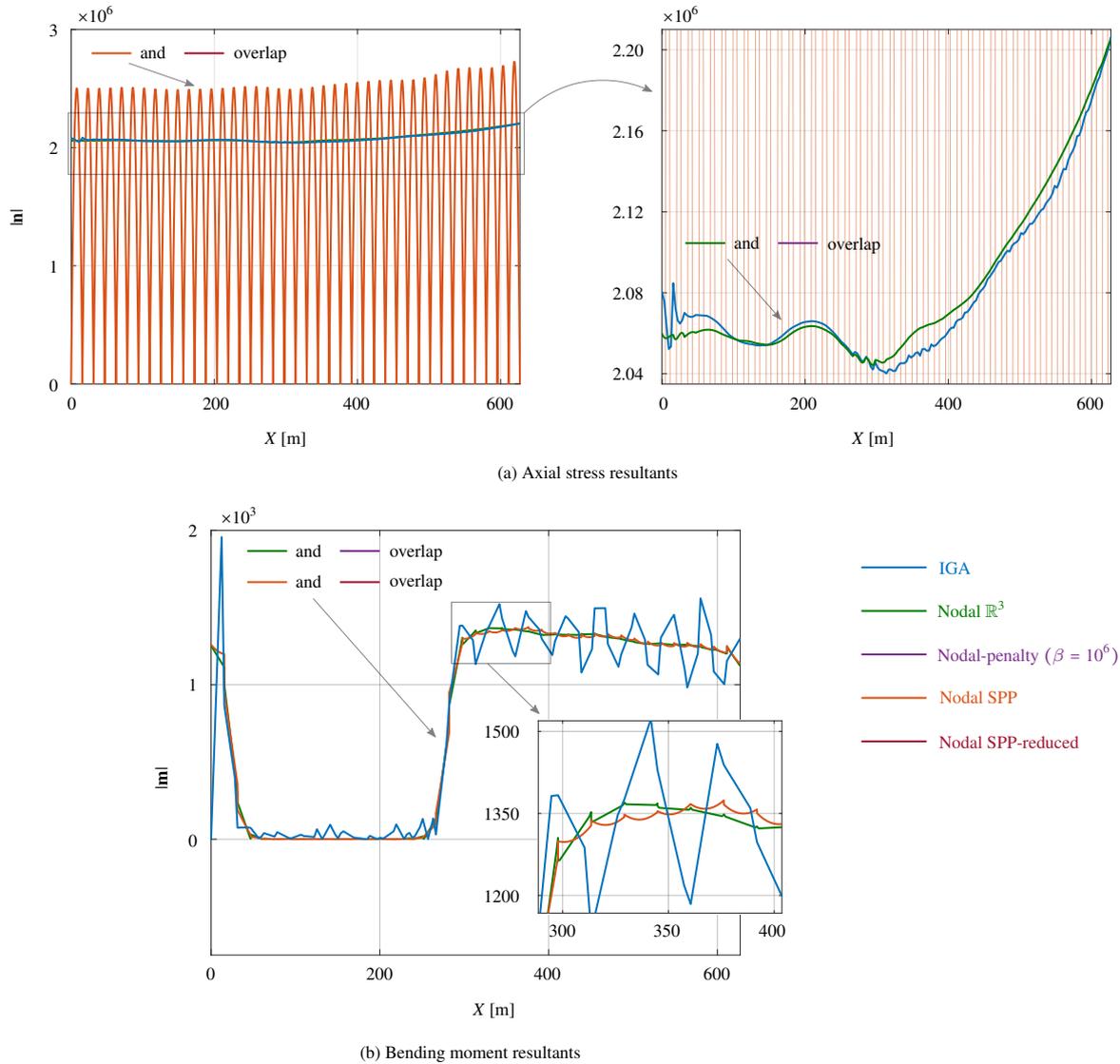    
	
	\caption{Axial stress and bending moment resultants in the final configuration (at $t=T$) of the mooring line in Figure \ref{fig:final_config_compare_mooringline}, computed with isogeometric and nodal discretization schemes.}\label{fig:stresses_compare_mooringline}
\end{figure}

\begin{figure}[ht]
	\centering
    \def\svgwidth{1\textwidth}
    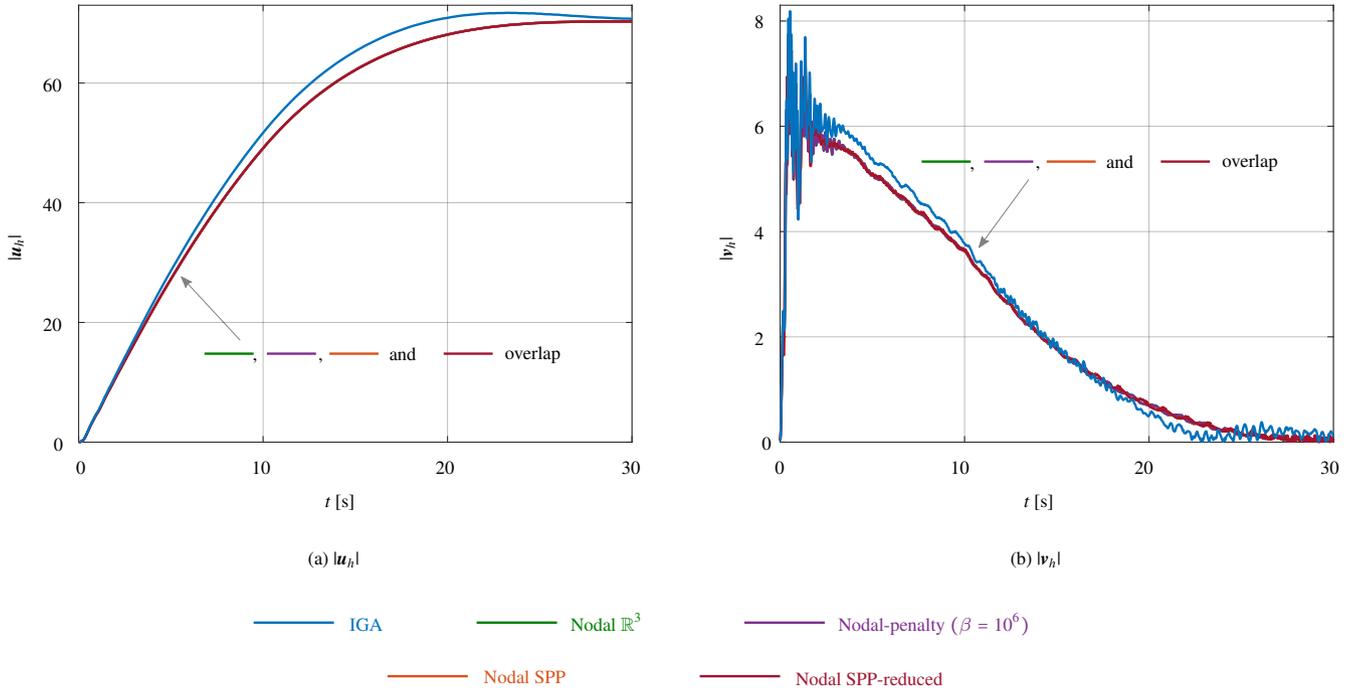    
	
	\caption{Time history of the displacement and velocity of the fairlead of the mooring line in Figure \ref{fig:final_config_compare_mooringline}, computed with isogeometric and nodal discretization schemes.}\label{fig:fairlead_compare_mooringline}
\end{figure}

\begin{figure}[ht]
	\centering
    \def\svgwidth{0.92\textwidth}
    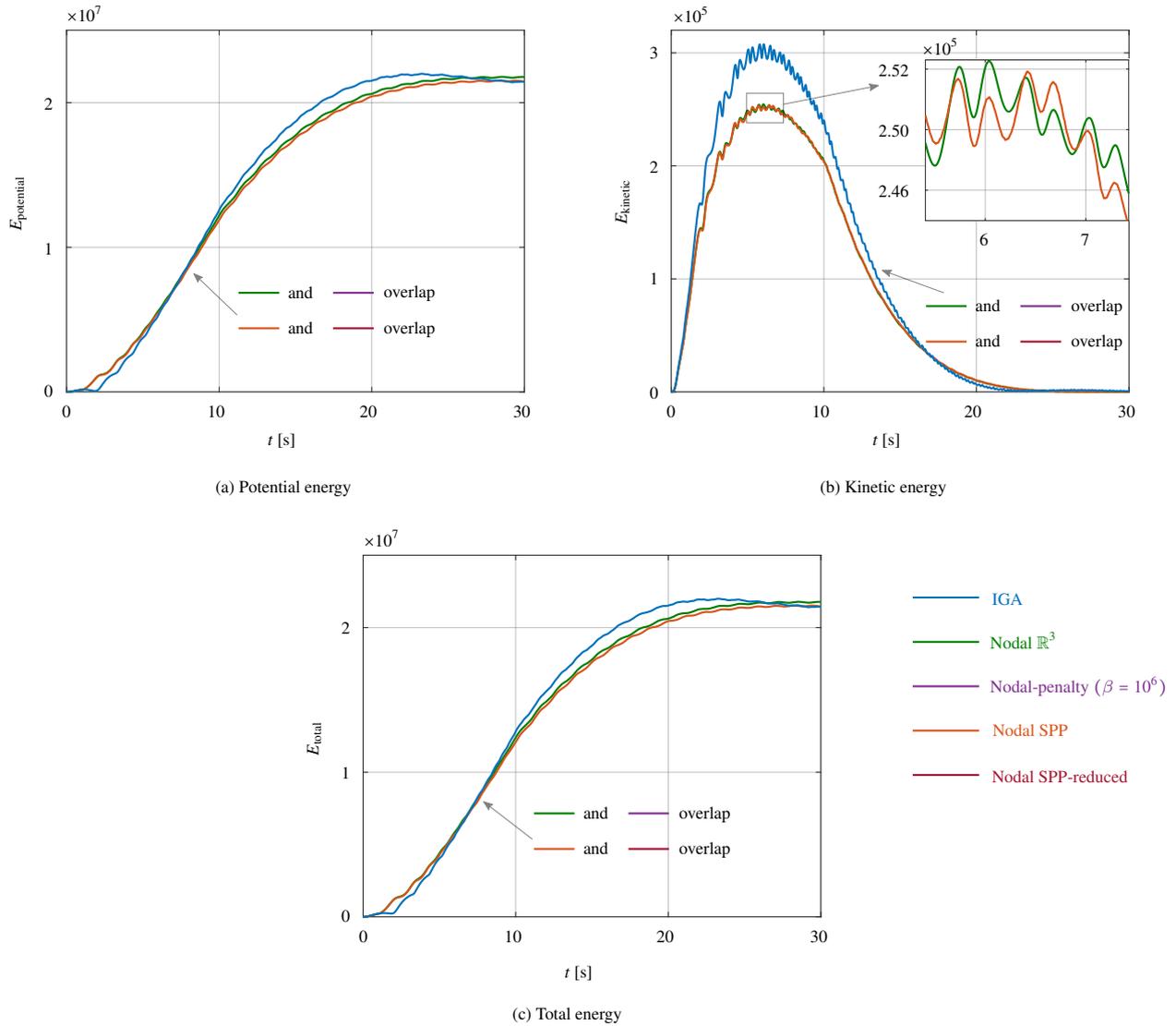    
	
	\caption{Time history of the energy of the mooring line in Figure \ref{fig:final_config_compare_mooringline}, computed with isogeometric and nodal discretization schemes.}\label{fig:energy_compare_mooringline}
\end{figure}

\begin{figure}[ht]
	\centering
    \def\svgwidth{0.85\textwidth}
    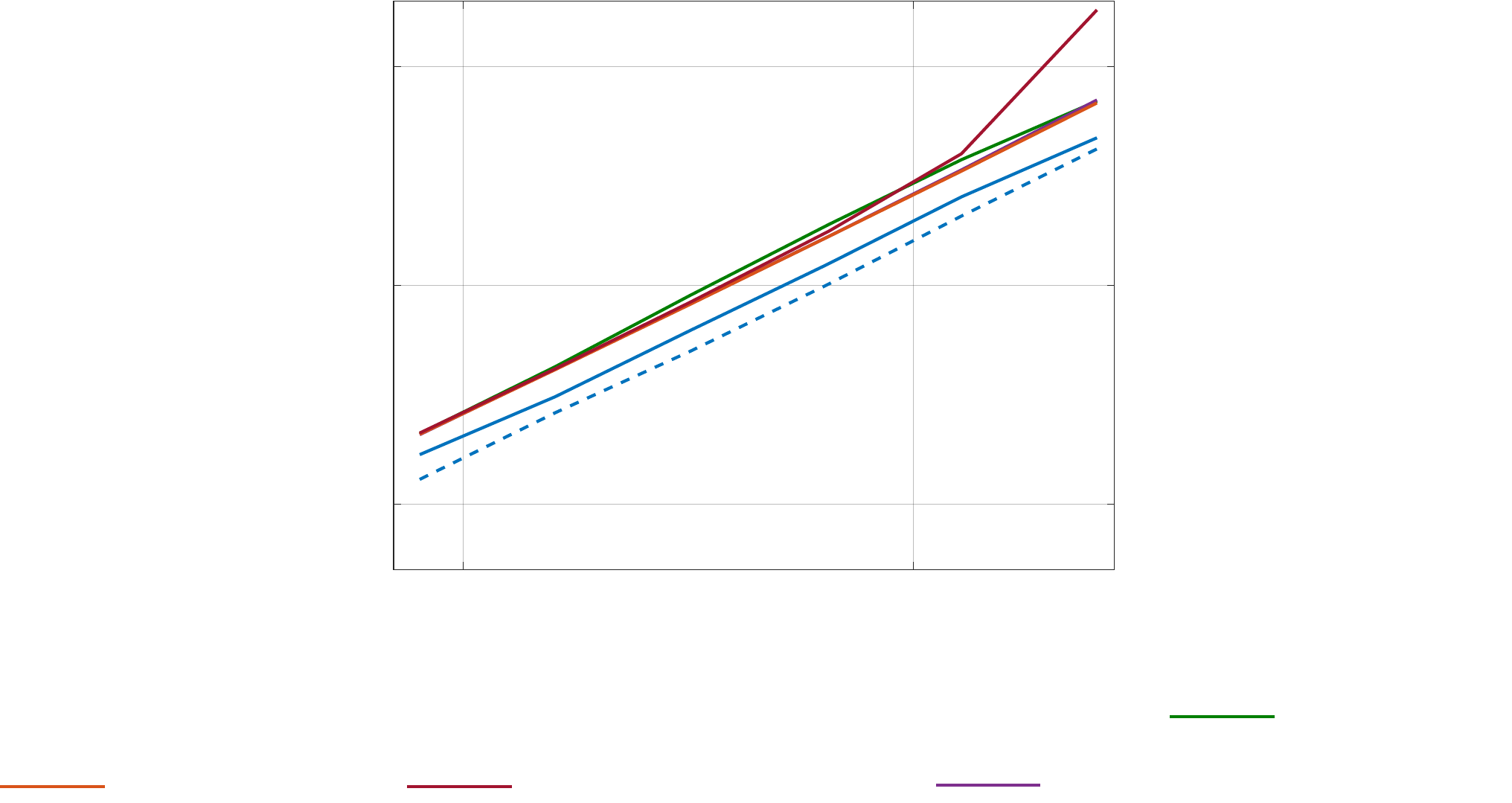    
	
	\caption{Computation time per iteration, averaged over all iterations and load steps, required when using isogeometric and nodal discretization schemes for the mooring line in Figure \ref{fig:final_config_compare_mooringline}.}\label{fig:time_compare_mooringline}
\end{figure}

\subsection{A dynamics example of mooring lines}\label{sec:results_dyn}

The second example of mooring lines considered in this work is a cable commonly employed in offshore wind engineering adapted from \cite[Ch.~7.6.4, p.~257]{Gunnar2024}. 
The cable has an 
initial length of $L_0=627$ m, 
a weight per unit length when submerged in water of $w_s=2.46$ kN/m, 
a Young modulus of $E = 5.6 \cdot 10^{10}$ N/m$^2$, 
and a cross-sectional area of $A = 0.0159$ m$^2$.  
The cable's left-end is fixed on the seabed while its right-end is also moved to a fairlead position with a prescribed point load in this case. 
This point load and the final configuration of the cable are provided by the author of \cite{Gunnar2024}.  
We note that F.G. Nielsen computed the final configuration by finding cable parts that rest on the seabed and deform due to the point load and cable's self-weights when submerged in water, considering the cable a sort of an elastic catenary and neglecting the effects of the water current. 
The provided solution is not based on a dynamic analysis but on a static one. 
In our computations, 
we want to capture and investigate the responses when the right-end is moved and the possible effects of the water current. 
Hence, we simulate the moving procedure as a dynamic analysis of the studied cable 
and consider the following logarithmic current speed profile (see also Figure \ref{fig:final_config_compare_mooringline}):
\begin{align}
    v_{\text{wind}} (z) = 2 \, \log \left( 1 + \frac{9}{z_{\text{seabed}}} \, z \right) \; [\text{m/s}] \,,
\end{align}
which can be employed for open water (see e.g. \cite{Montserrat_currentProfile2011}). 
We simulate the seabed at $z_{\text{seabed}} = 100$ m 
as a numerical barrier 
using the so-called \textit{barrier function} \cite{fiacco1990,wriggers2006,Roccia_mooring_line2024}. 
The main idea is to add a penalty term of the barrier function to the weak form \eqref{w-eom}, which introduces increasing energy/work when the distance between the rod and the barrier decreases and vice versa, 
i.e. the energy required to keep a distance between all point of the current configuration and the barrier. 
We note that this necessarily means that there is a numerical gap between the discrete rod configuration and the barrier, which can be regulated with the penalty factor associated with the barrier term. 
For more details on the barrier function and the linearization of associated terms for the studied rod formulation, we refer to \cite{Roccia_mooring_line2024}. 
For our computation, we consider a numerical barrier at $z=-0.5$ m such that a part of final configuration rests at a height of $z=0$. 
We choose a reciprocal function as the barrier function and an associated penalty factor of 25. 
We enforce the cable's weight and the prescribed point load at the right-end in 10 seconds with a time step of $0.01$ s. 
We continue the simulation with the constant final values of these forces in 20 seconds with the same time step, such that a part of the cable deforms and rests ``on'' the numerical barrier and reaches the final configuration. 
Similar to the static analysis above, 
we compute the response of this cable using the five semi-discrete formulations: \iga, \nodal, \nodalPenalty \ with $\beta=10^5$, \nodalSaddle, and \nodalSaddleRed. 
We discretize the cable in 40 elements and employ 
a tolerance of $10^{-10}$ for the convergence of the Newton-Raphson method in all cases. 
We again focus on the performance of different discretization schemes and hence employ basis functions in the same function space, i.e. cubic $C^1$ and Hermite splines. 
We enforce the homogeneous boundary condition on the left boundary in the same way as described in the static analysis above. 
The number of iterations is 4 for all studied formulations at each time step.

Figure \ref{fig:final_config_compare_mooringline} illustrates six snapshots every 5 seconds of the simulation, computed with cubic $C^1$ isogeometric discretization and outlier removal. 
We also include the provided final configuration of \cite[Ch.~7.6.4, p.~257]{Gunnar2024} (black solid curve) and the numerical barrier (black dashed horizontal line). 
Also for this example, we obtain indistinguishable configurations when using \iga \ without outlier removal or the other four formulations based on the nodal discretization scheme. 
Hence, we also plot the results of one formulation here purely for illustration clarity. 
Since we obtain the same responses when using \iga \ with and without outlier removal, we refer to \iga \ for both of these cases in the following. 
Comparing the configuration during the last 10 seconds with the provided solution of \cite[Ch.~7.6.4, p.~257]{Gunnar2024}, 
we see that they approximately overlap each other. 
This necessarily means that the effect of the considered current speed on the final configuration is negligible. 
Figure \ref{fig:stresses_compare_mooringline} illustrates the axial stress and bending moment resultants in the final configuration of the studied cable, computed with 
the aforementioned formulations. 
We also include an overkill solution (gray curve) computed with \iga\ using quintic ($p=5$) $C^4$ splines and 1024 elements as a reference solution. 
Focusing on axial stress resultants (Figure \ref{fig:stresses_compare_mooringline}a), 
we observe that 
using \iga \ (blue curve) or \nodal \ (green curve) leads to similar axial stresses that consist of oscillations, which might be caused by membrane locking and can be reduced by refining the mesh (see the gray curve in Figure \ref{fig:stresses_compare_mooringline}a). 
We see that using \nodalPenalty \ (purple curve) with a penalty factor of $\beta=10^5$ does not sufficiently enforce the unit nodal director constraint and thus leads to the same axial stress resultants as \nodal. 
Enforcing this constraint using either \nodalSaddle \ (orange curve) or \nodalSaddleRed \ (dark red curve), however, leads to zero nodal axial stresses, as observed and discussed in previous sections and examples. 
We note that also for this example, reducing the system of equations of \nodalSaddle \ using \nodalSaddleRed \ leads to the same stress resultants.

Focusing on the bending moment (Figure \ref{fig:stresses_compare_mooringline}b), we observe that using \iga \ leads to results consisting of oscillations with larger amplitude than that obtained with the other four formulations based on the nodal scheme. 
Using splines of higher polynomial and continuity order and refining the mesh reduces these oscillations only in parts of the mooring line with $X < 300$ m (see the gray curve in Figure \ref{fig:stresses_compare_mooringline}b). 
This necessarily means that membrane locking might not be the single cause of these oscillations in the bending moment resultants. 
We note that using any of the four formulations based on the nodal discretization scheme leads to bending moments that also oscillate but with significantly smaller amplitude. 
This might be caused by the lower accuracy in deformations in different error norms 
achieved with the chosen cubic $C^1$ splines and mesh size for this example, which was the case of the planar roll-up studied in Section \ref{sec:convergence_study_rollup}. 
The nodal values obtained with these formulations differ from each other due to the enforced unit nodal director constraint.

Figure \ref{fig:fairlead_compare_mooringline} illustrates the time history of the displacement and velocity at the right-end of the cable, i.e. the fairlead. 
We observe that all four formulations based on the nodal scheme lead to the same responses. 
Using \iga \ leads to slightly larger displacements between $t=10$ and $t \approx 28$ seconds and oscillating velocity with larger amplitude. 
Using splines of higher polynomial and continuity order and refining the mesh reduces these oscillations in the velocity (see the gray curve in Figure \ref{fig:fairlead_compare_mooringline}b). 
We see that at $t=30$ seconds, i.e. towards the end of the simulation, we obtain again approximately the same responses. 
This necessarily means that for this example, using \iga \ requires longer computation to reach the same final responses. 
We note that the reason for this observed difference might be the remaining outliers and/or activated high-frequency modes when using \iga \ (see also \cite{nguyen_rod_2024}). 
In Figure \ref{fig:energy_compare_mooringline}, we plot the time history of the potential, kinetic, and total energy obtained with the studied formulations. 
We observe that using the formulations based on the nodal scheme leads to approximately the same energy responses. 
We also see a slight difference in these responses when enforcing the unit nodal director constraint using the Lagrange multiplier method. 
This is expected 
since this constraint leads to different nodal stress resultants. 
Focusing on the blue curves obtained with \iga, 
we observe that 
it leads to slightly larger energy than using the nodal scheme, which is due to slightly larger responses observed in Figure \ref{fig:fairlead_compare_mooringline}. 
Using splines of higher polynomial and continuity order and refining the mesh leads to the same energy responses (see the gray curve in Figure \ref{fig:energy_compare_mooringline}). 
Towards the end of the computation, we again obtain approximately the same energy responses as using the nodal discretization scheme.

In terms of computational cost, we also investigate the number of iterations and computing time per iteration when using the studied formulations under mesh refinement. 
For this dynamic example, all formulations require the same number of iterations per time step that is 4 iterations per time step. 
In 
Figure \ref{fig:time_compare_mooringline}, we plot the averaged computing time per iteration and time step when using the studied formulations. 
We have the same observations as in the static analysis above: 
\iga \ (blue curves) requires the least computing time per iteration despite the global matrix multiplication for enforcing the boundary conditions and employing the outlier removal approach \cite{hiemstra_outlier_2021}. 
Using \nodal, \nodalSaddle, \nodalPenalty, and \nodalSaddleRed \ requires the same computing time per iteration. 
One exception is \nodalSaddleRed \ for computations on fine meshes, which 
requires significantly more computational effort due to the reassembly of the nullspace matrix per iteration.

We conclude that for the dynamic analysis of the exemplary mooring line, all four formulations based on the nodal discretization scheme lead to approximately the same responses, except the axial stress resultants which are constrained to zero nodal value when enforcing the unit nodal director constraint using the Lagrange multiplier method. 
For this example, using cubic $C^1$ isogeometric discretization leads to approximately the same final configuration, however, different stress resultants, particularly bending moments with larger oscillations. 
It also shows slightly larger displacement and velocity responses during parts of the simulation, which may be due to the remaining outliers and/or activated high-frequency modes. 
Hence, using cubic $C^1$ \iga \ may require longer computation, finer meshes, or cubic splines of higher continuity order to obtain the same responses as using the nodal scheme. 
Nevertheless, \iga \ requires the least time per iteration, with or without outlier removal, despite the global matrix multiplication. 
Using any of the four formulations based on the nodal scheme requires approximately the same computing time, except \nodalSaddleRed \ which requires significantly more time on finer meshes.

\section{Summary and conclusions}\label{sec:summary}
  
In this work, 
we discussed and compared the nodal and isogeometric spatial discretization schemes for the nonlinear formulation of shear- and torsion-free rods \cite{gebhardt_2021_beam}. 
We showed that while the latter leads to a discrete solution in multiple copies of $\mathbb{R}^3$, the former leads to a discrete solution in multiple copies of either the same space or the manifold $\mathbb{R}^3 \times S^2$. 
Preserving the unit sphere $S^2$ structure of the director field at the nodes leads to a discrete solution in multiple copies of $\mathbb{R}^3 \times S^2$ and requires an additional unit nodal director constraint, 
which leads to zero nodal axial stress values, i.e. oscillating axial stress resultants. 
We studied five semi-discrete formulations and corresponding matrix equations of different discretization variants: 
isogeometric discretizations (\iga), 
nodal discretization without considering unit nodal director constraint (\nodal), 
nodal discretization with a nodal director constraint enforced using Lagrange multiplier method (\nodalSaddle), 
nodal discretization with a reduced system of equations of \nodalSaddle \ using nullspace method (\nodalSaddleRed), and
nodal discretization with a nodal director constraint enforced using the penalty method (\nodalPenalty). 
We discussed the computational cost related to each of these five formulations and showed that  
\nodalSaddle \ leads to the largest system of equations in the form of a saddle-point problem, compared to the other three formulations based on the nodal scheme. 
\iga \ enables quadratic and higher-order continuous basis functions, possibly leading to a smaller system. 
We also showed that 
formulations with a unit nodal director constraint enforced using the Lagrange multiplier method lead to a non-symmetric system matrix for dynamic computations, possibly requiring more computational effort for solving significantly large system of equations. 
We numerically illustrated via examples of a catenary and mooring line that all studied formulations generally require the same number of iterations and using \iga \ requires the least time per iteration, with or without outlier removal. 
Using formulations based on the nodal scheme requires approximately the same computing time, except \nodalSaddleRed \ which requires significantly more time on fine meshes due to the reassembly of the nullspace matrix in each iteration.

Furthermore, we numerically studied the condition number of the resulting system matrix, gaining insights into the robustness of the studied formulations. 
We illustrated for an exemplary cable that \iga \ leads to slightly smaller condition number than \nodal, which is similar to that obtained with \nodalSaddleRed. 
We showed for this example that employing the scaled director conditioning \cite{Kloppel_scaled_director2011,Wall_scaled_director2000,Gee_scaled_director2005} improves the conditioning and hence the robustness when using \nodalSaddle \ or \nodalPenalty. 
We then numerically illustrated via a pure bending example of a planar roll-up that preserving the nodal director field in the unit sphere leads to better accuracy in the deformations in different error norms. 
Our numerical results imply the effect of membrane locking on all studied semi-discrete formulations, particularly oscillations in stress resultants and convergence of errors in the deformations. 
We showed for two examples of a catenary and a mooring line that all formulations approximately lead to the same final deformed configuration. 
For the dynamics example, however, 
cubic $C^1$ isogeometric discretization leads to slightly larger responses and hence 
may require longer computation to reach the same final responses as using the nodal discretization scheme. 
This may be due to the remaining outliers and/or activated high-frequency modes.

Based on the presented results, nodal discretizations may be preferable for dynamic computations due to their higher robustness compared to cubic isogeometric discretizations. 
Comparing the four studied semi-discrete formulations when using nodal discretizations, \nodalSaddleRed\ showed greater robustness than the others and required approximately the same computation time per time step when using coarser meshes. 
When discretizing with very fine meshes, \nodalPenalty\ may be a preferable alternative that is also robust and requires less computational effort. 
For static analyses, our results indicated that not only \nodalSaddleRed\ and \nodalPenalty, but also \iga, are robust and require approximately the same number of iterations. However, \nodalPenalty\ and \iga\ generally require less computation time per iteration and may therefore be more preferable. 
When choosing a discretization scheme for either static or dynamic computations, the accuracy of the responses, particularly the stress resultants, also plays a crucial role. 
Based on the presented results, \iga\ and \nodalPenalty\ may be preferable for static and dynamic computations, respectively, as they offer a favorable compromise between unconstrained or weakly constrained nodal axial stress values and robustness compared to the other studied formulations.

The results presented in this work give an overview and deeper understanding of different spatial discretization schemes for the nonlinear rod formulation \cite{gebhardt_2021_beam}. There are a number of avenues for future work. 
One aspect is to investigate and eliminate the effect of membrane locking on the studied semi-discrete rod formulations. 
To this end, we are particularly interested in approaches such as 
reduced/selective integration \cite{Noor1981,Adam2015,Zou_quadrature2021} or 
those based on variational multiscale stabilization (see e.g. \cite{Aguirre_VMS2023}). 
For the considered nonlinear rod formulation, it is particularly interesting to identify whether other locking phenomena occur, such as those reported in \cite{Willmann_nonlinear_locking2023} for geometrically nonlinear shell structures.  
A second aspect for future work is the development of other strain measures that can address zero nodal axial stress values while preserving the nodal director field in the unit sphere for nodal discretizations. 
Another aspect is to develop and investigate an analytical scheme to estimate the scaling penalty factor when using nodal discretizations and the penalty method to enforce the unit nodal director constraint. Such a scheme allows the chosen penalty factor to be a problem-independent intensity factor\cite{Pasch2021}.

\bmsection*{Acknowledgments}

T.-H. Nguyen, B.A. Roccia, and C.G. Gebhardt gratefully acknowledge the financial support from the European Research Council through
the ERC Consolidator Grant “DATA-DRIVEN OFFSHORE” (Project ID 101083157). 
D. Schillinger 
gratefully acknowledges financial support from the German Research Foundation (Deutsche Forschungsgemeinschaft)
through the DFG Emmy Noether Grant SCH 1249/2-1 and the standard DFG grant SCH 1249/5-1.

\appendix

\section{Variation and linearization of the unit nodal director constraint}\label{sec:linearized_constraint_vector}

We recall Equation \eqref{eq:variation_of_constraint} that defines the matrix $\mat{J}$ resulting from the variation of $\vect{\Psi}$:
\begin{align}
    & \delta \boldsymbol{\Psi} = 2\,\left[\delta \vect{d}_1 \cdot \vect{d}_1 \; \ldots \; \delta \vect{d}_{n_n} \cdot \vect{d}_{n_n} \right]^T = \mat{J} \, \delta \, \Bar{\vect{\qhat}}\,,
\end{align}
where $n_n$ is the number of nodes. 
This necessary means: 
\begin{align}
    & \mat{J} = 2 \, \begin{bmatrix}
        \mat{0} & \mat{d}_1^T & \mat{0} & \ldots &  \\
         &  & \mat{0} & \mat{d}_2^T & \mat{0} & \ldots  \\
         &  &  & \ldots & & \\
         &  &  &  & \mat{0} & \mat{d}_{n_n}^T
    \end{bmatrix}
    \,. 
\end{align}
The nullspace matrix $\mat{D}$ of $\mat{J}$, such that $\mat{J} \, \mat{D} = \mat{0}$, is:
\begin{align}\label{eq:nullspace}
    \mat{D} =
    \begin{bmatrix}
        \vect{E}_1 & \vect{E}_2 & \vect{E}_3 & \mat{0} & \mat{0} & \ldots & & & \\
         &  & \mat{0} & \hat{\mat{d}}_1^2 & \hat{\mat{d}}_1^3 & \mat{0} & \ldots & & \\
         &  & & & \mat{0} & \vect{E}_1 & \vect{E}_2 & \vect{E}_3 & \ldots \\
         &  & & &  &  & & \ldots & \ldots & \\
         &  & & &  &  & & \mat{0} & \hat{\mat{d}}_n^2 & \hat{\mat{d}}_n^3 
    \end{bmatrix}
    \,, 
\end{align}
where $\hat{\mat{d}}_i^2$ and $\hat{\mat{d}}_i^3$ are two dual vectors of the $i$-th nodal director $\mat{d}_i$ and are computed as: $\hat{\mat{d}}_i^j = [\mat{d}_i]_\times \, \vect{E}_j$, $j = 1,2,3$, i.e.
\begin{align}
    \hat{\mat{d}}_i^1 = \begin{bmatrix}
        0 \\ d_i^3 \\ - d_i^2
    \end{bmatrix} \,, \quad \hat{\mat{d}}_i^2 = \begin{bmatrix}
        - d_i^3 \\ 0 \\ d_i^1
    \end{bmatrix} \,, \quad \hat{\mat{d}}_i^3 = \begin{bmatrix}
        d_i^2 \\ - d_i^1 \\ 0
    \end{bmatrix} \,.
\end{align}
Here, 
$d_i^j$ is the $j$-th component of the director $\mat{d}_i$. 
We note that $\hat{\mat{d}}_i^j$, $j = 1,2,3$, are linear dependent. In particular, one of these three dual vectors can always be expressed as a linear combination of the other two. 
Hence, an arbitrary pair of these three vectors consists two linear independent vectors. 
For our computations in this work, we choose $\hat{\mat{d}}_i^2$ and $\hat{\mat{d}}_i^3$ for $\mat{D}$.

Consider the matrix equations \eqref{eq:matrix_eq_nodalR3S2-strong}. 
The counterpart $\mat{A}_c = \mat{A}_c(\vect{\lambda})$, i.e. the contribution of the unit nodal director constraint to the system matrix, 
results from the linearization of the term $\delta \Bar{\vect{\qhat}} \cdot \mat{J}^T (\Bar{\vect{\qhat}}) \,\vect{\lambda}$ evaluated at $t_{n+\frac{1}{2}}$, but with respect to $\Bar{\vect{\qhat}}_{n+1}$. 
$\mat{A}_c$ is then: 
\begin{align}
    & \mat{A}_c = \mat{A}_c\left(\vect{\lambda}_{n+\frac{1}{2}}\right) = \begin{bmatrix}
        \mat{0} & \lambda_1 \mat{I} & & &  \\
         & \mat{0} & \lambda_2 \mat{I} & & \\
         & & & \ldots & \\
         & & &\mat{0} & \lambda_{n_n} \mat{I}
    \end{bmatrix}
    \,.
\end{align}

Consider the matrix equations \eqref{eq:matrix_eq_nodalR3S2-strong-reduced}.
The counterpart $\mat{A}_D$, i.e. the contribution of the nullspace matrix to the system matrix, 
results from the linearization of the nullspace matrix $\mat{D}^T_{n+\frac{1}{2}}$ with respect to $\Bar{\vect{\qhat}}_{n+1}$. 
The $i$-th column of $\mat{A}_D$, $i=1,\ldots,6n_n$, is then:
\begin{align}
    \mat{A}_{D,i} = \left(\Delta_i \, \mat{D}_{n+\frac{1}{2}}\right)^T \left(\Bar{\vect{F}}^{\text{ext}}_{n+\frac{1}{2}} - \Bar{\vect{F}}\right) \,,
\end{align}
where $\Delta_i \, \mat{D}_{n+\frac{1}{2}}$ is the linearization of $\mat{D}_{n+\frac{1}{2}}$ with respect to the $i$-th degree of freedom (dof). 
Since $\mat{D}$ only depends on the nodal director, not the nodal position (see \eqref{eq:nullspace}), 
one needs to compute $\mat{A}_D$ for only three dofs at each $j$-th node, $j=1,\ldots,n_n$, that are: $i = 6j-2,\, 6j-1, \, 6j$. 
For these dofs, the non-zero $3\times2$ block matrix of $\Delta_i \, \mat{D}_{n+\frac{1}{2}}$ expands on the $(6j-2:6j)$-th rows and $(5j-1:5j)$ columns and takes the following form:
\begin{subequations}
\begin{align}
    & \Delta_{6j-2} \mat{D}= \frac{1}{2} \, \begin{bmatrix}        
        \vect{E}_3 & -\vect{E}_2
    \end{bmatrix} \,, \\
    & \Delta_{6j-1} \mat{D} = \frac{1}{2} \, \begin{bmatrix}        
        \mat{0} & \vect{E}_1
    \end{bmatrix} \,,  \\
    & \Delta_{6j} \mat{D} = \frac{1}{2} \, \begin{bmatrix}        
        -\vect{E}_1 & \mat{0}
    \end{bmatrix} \,. 
\end{align}    
\end{subequations}
Here, $\vect{E}_i$, $i=1,2,3$, are the canonical Cartesian basis of $\mathbb{R}^3$. 
We note that the factor $\frac{1}{2}$ results from the chain rule employed when linearizing $\mat{D}_{n+\frac{1}{2}}$ with respect to $\Bar{\vect{\qhat}}_{n+1}$.

Consider the matrix equations \eqref{eq:matrix_eq_nodalR3S2-weak}.
The counterpart $\mat{A}_{\beta}$, i.e. the contribution of  
the penalty term to the system matrix, 
results from 
the linearization of the penalty term $\left( \beta \frac{2EI}{L} \, \delta \Bar{\vect{\qhat}} \cdot \mat{J}^T (\Bar{\vect{\qhat}}) \,\vect{\Psi} \right)$ evaluated at the time instance $t_{n+1}$ with respect to $\Bar{\vect{\qhat}}_{n+1}$. 
$\mat{A}_\beta$ is then: 
\begin{align}
    \mat{A}_\beta = \mat{A}_\beta \left(\Bar{\vect{\qhat}}_{n+1}\right) & = \beta \frac{2EI}{L} \, \mat{J}^T \, \mat{J} \nonumber \\
    & +\, 2\beta \frac{2EI}{L}\, 
    \begin{bmatrix}
        \mat{0} & \left(\vect{d}_1 \cdot \vect{d}_1 - 1 \right) \mat{I} & & &  \\
         & \mat{0} & \left(\vect{d}_2 \cdot \vect{d}_2 - 1 \right) \mat{I} & & \\
         & & & \ldots & \\
         & & &\mat{0} & \left(\vect{d}_{n_n} \cdot \vect{d}_{n_n} - 1 \right) \mat{I}
    \end{bmatrix} 
    \,.
\end{align}

\section{Convergence study of a planar roll-up with different slenderness ratios}\label{sec:convergence_rollup_radi}

\begin{figure}[h]
	\centering
    \def\svgwidth{0.831\textwidth}
    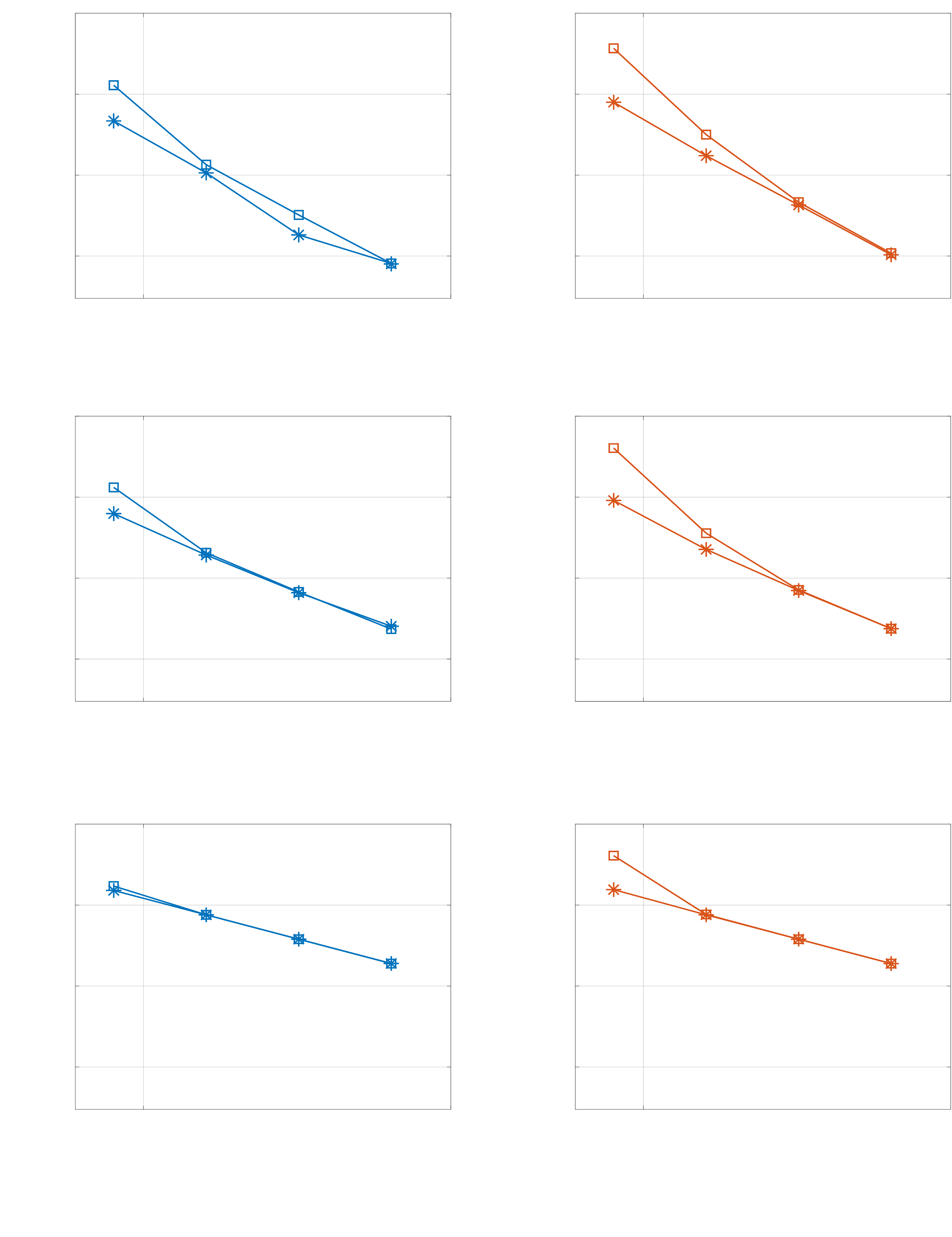

 \caption{Convergence of relative errors of the clamped rod bent to a circle computed with \textbf{cubic $C^1$} isogeometric discretizations and \nodalSaddle.} \label{fig:convergence_compare_rollup_radi1}
\end{figure}

\begin{figure}[h]
	\centering
    \def\svgwidth{0.95\textwidth}
    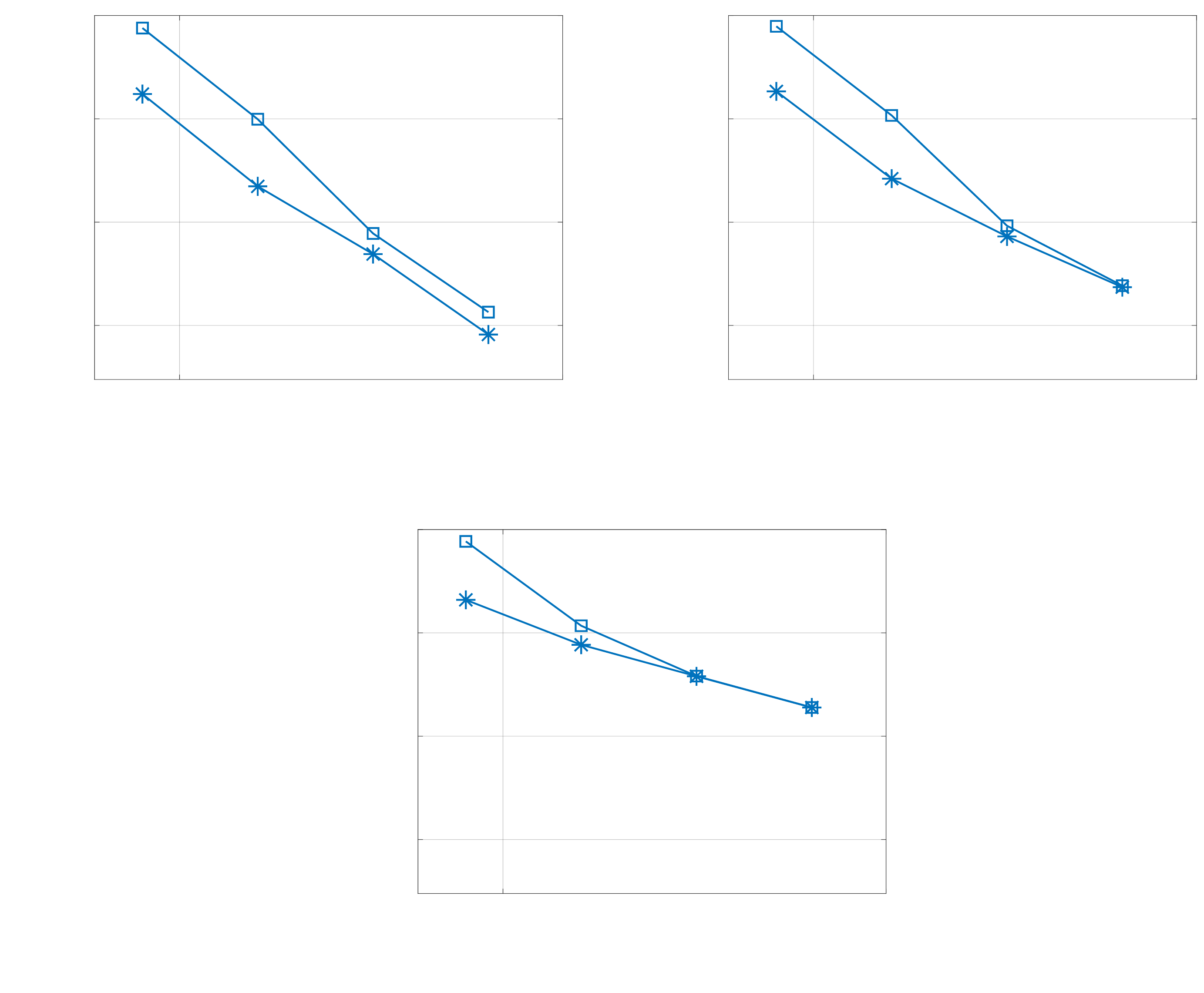

 \caption{Convergence of relative errors of the clamped rod bent to a circle computed with \textbf{cubic $C^2$} isogeometric discretizations.} \label{fig:convergence_compare_rollup_radi2}
\end{figure}

For the pure-bending example of a planar roll-up studied in Section \ref{sec:convergence_study_rollup}, 
we now investigate the accuracy and convergence behavior obtained with the studied semi-discrete formulations for different slenderness ratios. 
We consider the planar roll-up with a circular cross-section and hence employ the ratio of the initial length $L$ to the cross-sectional radius $R$, i.e. a slenderness ratio $L/R$. 
In Figure \ref{fig:convergence_compare_rollup_radi1}, we plot the convergence of the relative errors in $L^2$-, $H^1$-, and $H^2$-norms of the studied planar roll-up, obtained with cubic $C^1$ isogeometric discretization and \nodalSaddle. 
To investigate possible effect of higher-order continuous splines in this context, we also plot the convergence curves obtained with cubic $C^2$ isogeometric discretization in Figure \ref{fig:convergence_compare_rollup_radi2}. 
We note that here, we also obtain the same accuracy and convergence behavior in all three error norms when using \nodalPenalty, \nodalSaddle, and \nodalSaddleRed, while using \nodal \ again leads to ill-conditioned system matrix for certain meshes. 
Hence, we illustrate only the results obtained with \nodalSaddle \ (Figures \ref{fig:convergence_compare_rollup_radi1}b, d, and f) and compare with those obtained with \iga. 
We observe a pre-asymptotic plateau in the convergence curves with large slenderness ratios, which becomes more severe with increasing slenderness ratio, in all error norms when using either isogeometric or nodal discretization scheme. 
This result implies the well-known effect of membrane locking on the accuracy and convergence behavior of employed discretizations. 
Comparing the error levels in all three error norms obtained with cubic $C^1$ isogeometric discretization and \nodalSaddle \ (see Figure \ref{fig:convergence_compare_rollup_radi1}), we see that membrane locking leads to higher error levels when using \nodalSaddle, i.e. its effect on \nodalSaddle \ is more severe than on cubic $C^1$ \iga. 
Using cubic $C^2$ \iga, however, leads to the same error levels as \nodalSaddle \ (see Figure \ref{fig:convergence_compare_rollup_radi2}). 
We note that there are various locking-preventing approaches, that are well-established for both isogeometric and nodal discretization, 
for example, the approach of reduced/selective integration \cite{Noor1981,Adam2015,Zou_quadrature2021} or approaches based on Hu-Washizu or Hellinger–Reissner variational principles \cite{Cannarozzi2008,Choit1995,Lee1993}. 
We plan to investigate the locking effect and locking-preventing techniques for the studied semi-discrete rod formulations in future work.

\section{A dynamics three-dimensional example of mooring lines}\label{sec:3Ddynamic_example}

\begin{figure}[ht]
	\centering
    \def\svgwidth{0.7\textwidth}
    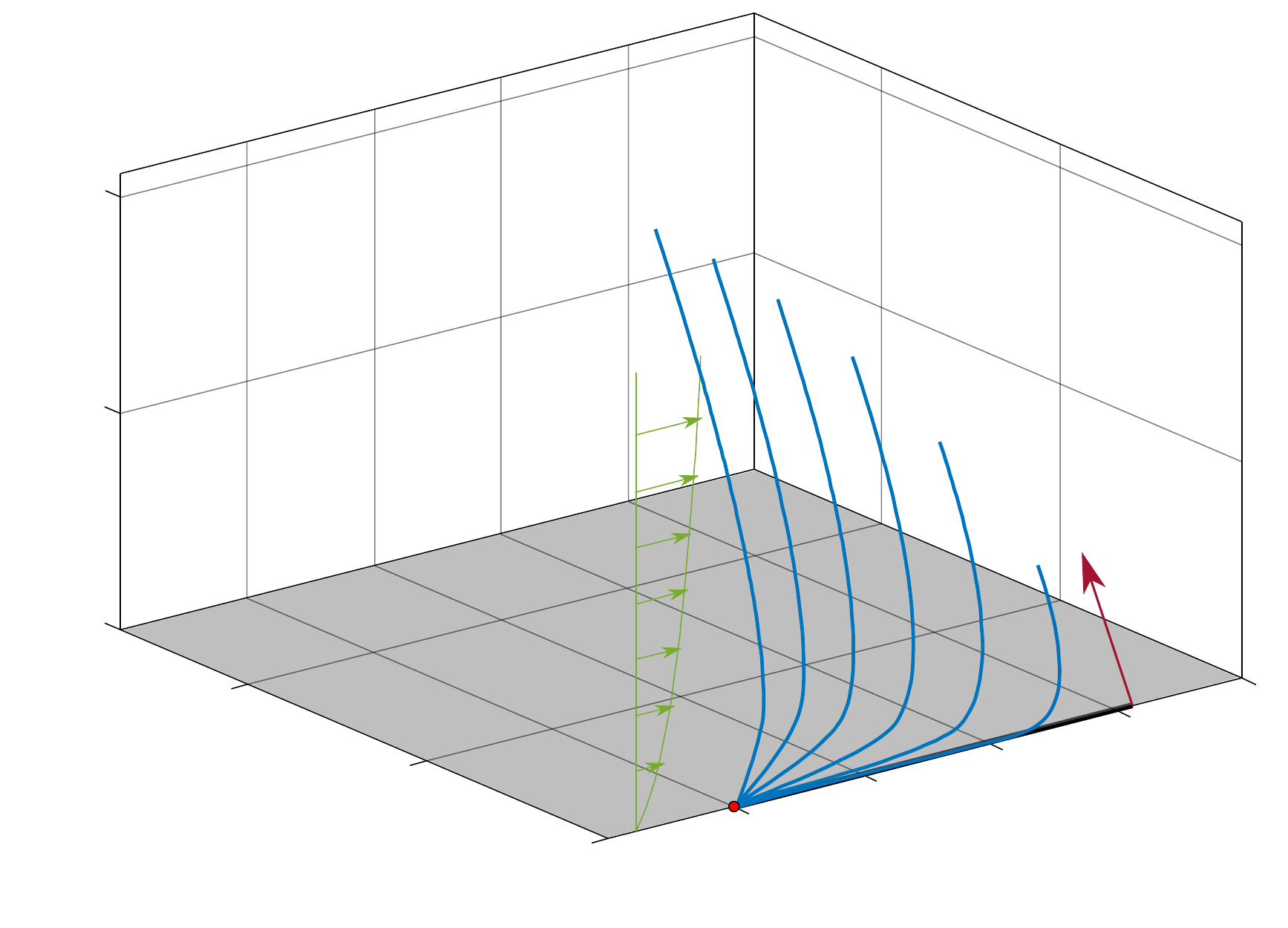

 \caption{Snapshots of a mooring line subjected to a logarithmic current speed profile in open water and a single force at the fairlead, computed with cubic $C^1$ isogeometric discretization and the outlier removal approach \cite{hiemstra_outlier_2021}. The horizontal black dashed line represents the numerical barrier.} \label{fig:final_config3D_compare_mooringline}
\end{figure}

\begin{figure}[ht]
	\centering
    \def\svgwidth{0.87\textwidth}
    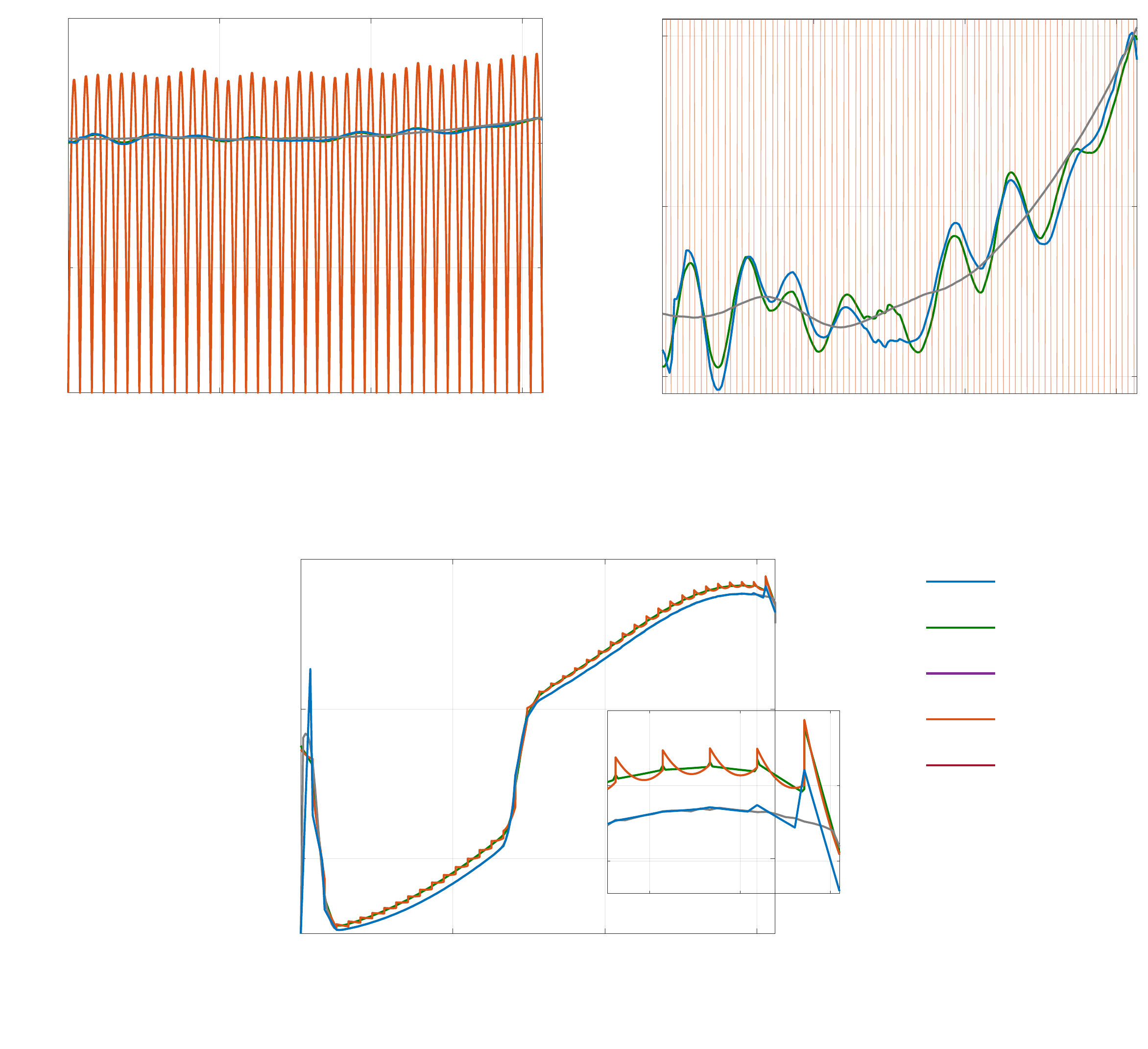    
	
	\caption{Axial stress and bending moment resultants in the final configuration (at $t=T$) of the mooring line in Figure \ref{fig:final_config3D_compare_mooringline}, computed with isogeometric and nodal discretization schemes.}\label{fig:stresses3D_compare_mooringline}
\end{figure}

\begin{figure}[ht]
	\centering
    \def\svgwidth{1\textwidth}
    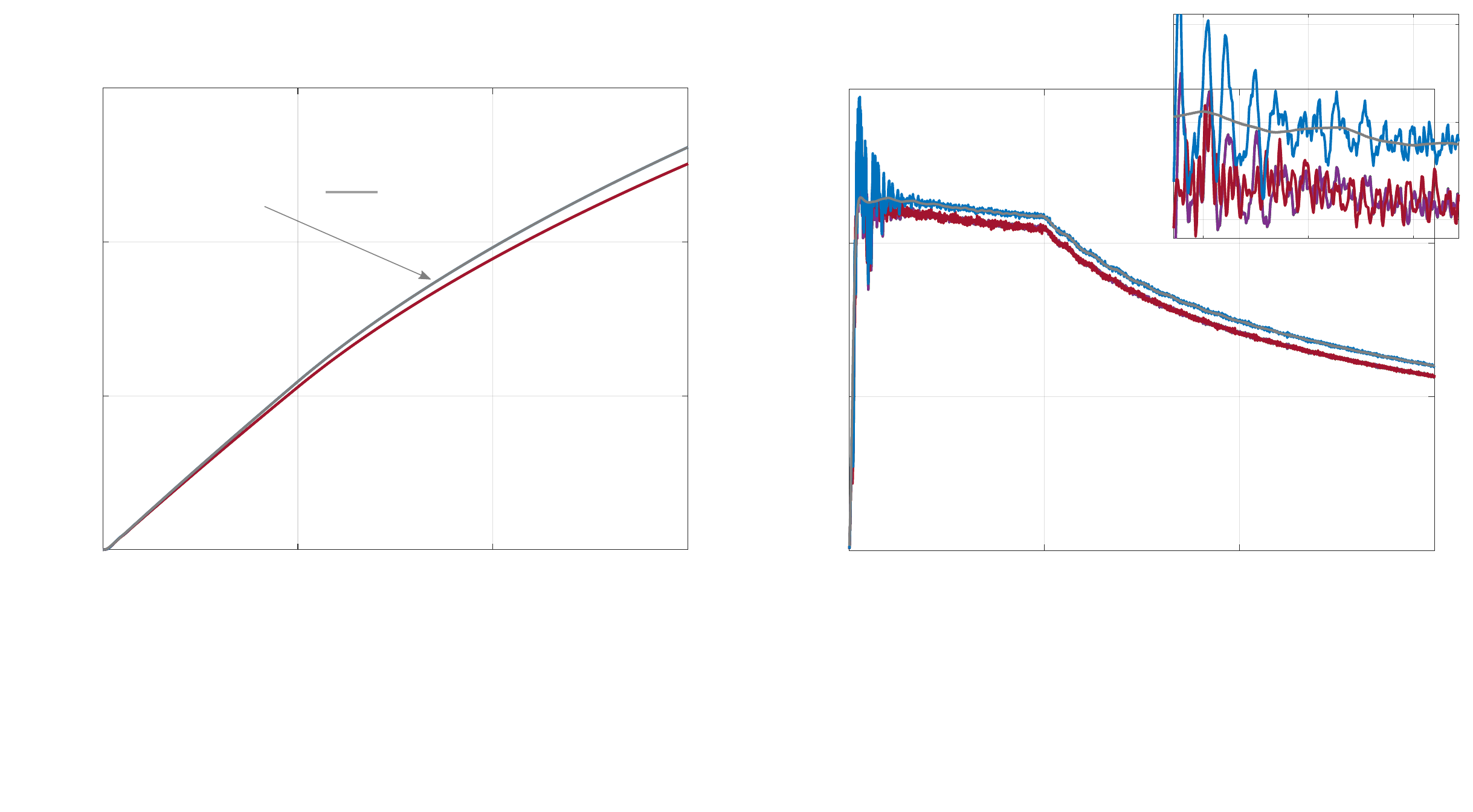    
	
	\caption{Time history of the displacement and velocity of the fairlead of the mooring line in Figure \ref{fig:final_config3D_compare_mooringline}, computed with isogeometric and nodal discretization schemes.}\label{fig:fairlead3D_compare_mooringline}
\end{figure}

\begin{figure}[ht]
	\centering
    \def\svgwidth{0.92\textwidth}
    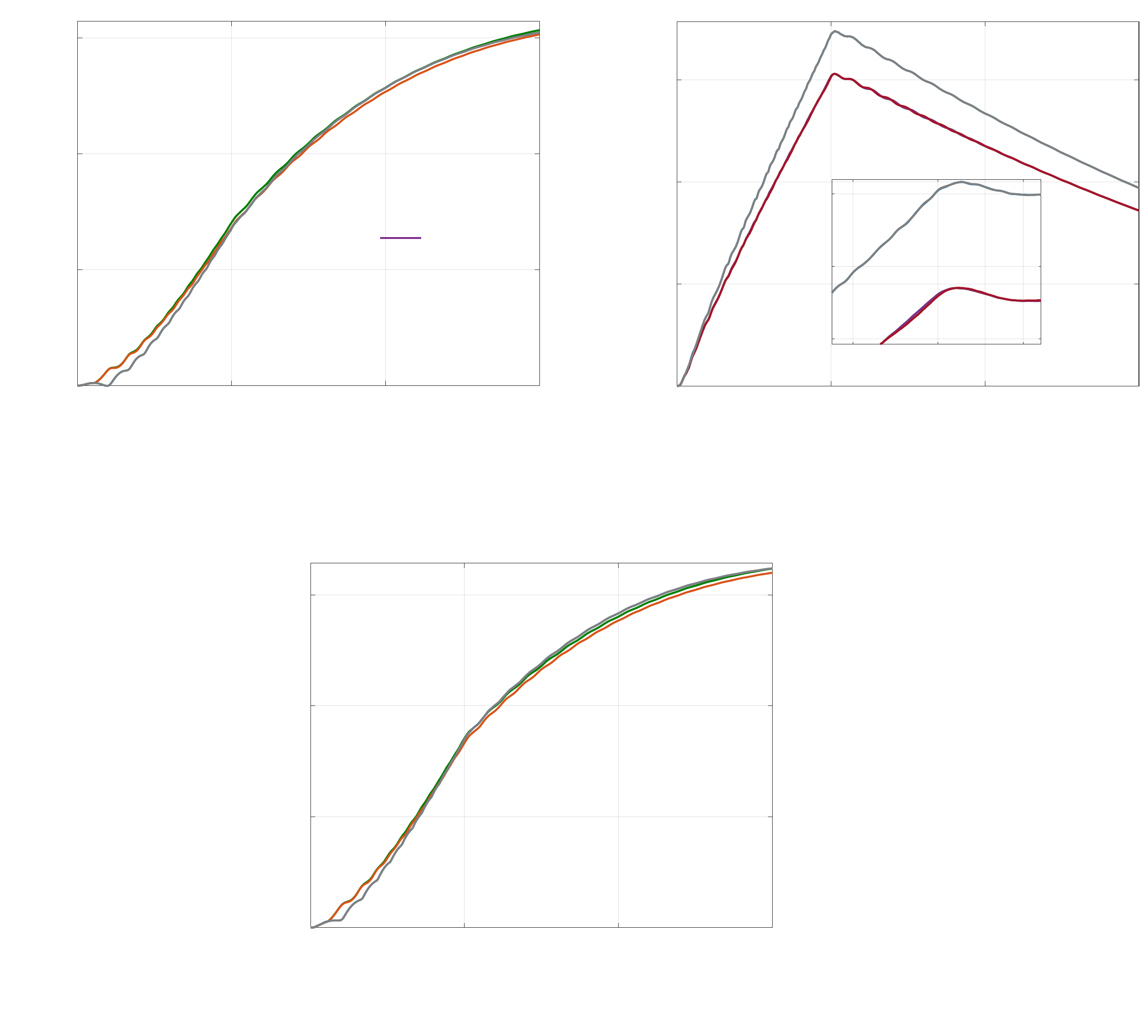    
	
	\caption{Time history of the energy of the mooring line in Figure \ref{fig:final_config3D_compare_mooringline}, computed with isogeometric and nodal discretization schemes.}\label{fig:energy3D_compare_mooringline}
\end{figure}

In this work, we also consider a three-dimensional example of mooring lines that is a modification of the example studied in Section \ref{sec:results_dyn}. 
Instead of an in-plane point load at the fairlead, $F_{end} = \left[F_x, 0, F_z \right]$, as considered in Section \ref{sec:results_dyn}, 
we apply a three-dimensional point load at the fairlead of the same magnitude, $\left[F_x/\sqrt{2}, F_x/\sqrt{2}, F_z \right]$, as illustrated in Figure \ref{fig:final_config3D_compare_mooringline}. 
We again employ the same geometry, cable material, and current profile as in Section \ref{sec:results_dyn}. 
We also consider the seabed as a numerical barrier and apply the same barrier function as in Section \ref{sec:results_dyn}, which naturally extend the barrier to a flat surface (see also Figure \ref{fig:final_config3D_compare_mooringline}). 
We compute the final configuration and responses in the same way as in Section \ref{sec:results_dyn}, i.e. a computation of 20 seconds with constant loading after enforcing the point load at the fairlead in a stepwise manner within 10 seconds. 
We also compare the responses obtained with the five semi-discrete formulations: \iga, \nodal, \nodalPenalty \ with $\beta=10^5$, \nodalSaddle, and \nodalSaddleRed, using the same discretizations as in Section \ref{sec:results_dyn}. The number of iterations is 4 for all studied formulations at each time step.

Figure \ref{fig:final_config3D_compare_mooringline} illustrates six snapshots every 5 seconds of the simulation, computed with cubic $C^1$ isogeometric discretization and outlier removal. 
For illustration clarity, we only plot the response obtained with \iga\ since that obtained with other four formulations are indistinguishable for this example. 
Figure \ref{fig:stresses3D_compare_mooringline} illustrates the axial stress and bending moment resultants in the final configuration of the studied cable, computed with 
the five aforementioned formulations. 
We include an overkill solution (gray curve) computed with \iga\ using quintic ($p=5$) $C^4$ splines and 1024 elements as a reference solution. 
We observe the same results as in the two-dimensional example studied in Section \ref{sec:results_dyn}. 
In 
Figure \ref{fig:fairlead3D_compare_mooringline} and 
\ref{fig:energy3D_compare_mooringline}, we plot the time history of the responses at the fairlead, 
displacement and velocity, and the energy 
(potential, kinetic, and total energy), respectively, obtained with the studied formulations. 
We again observe the same results as in the example studied in Section \ref{sec:results_dyn}.

\bibliography{wileyNJD-AMA.bib}

\end{document}